\documentclass{article}

\usepackage{amssymb,amsmath,theorem,euscript}

\input{epsf}

\newcounter{sec}

\newcounter{punct}[sec]

\def\punct{\refstepcounter{punct}{\arabic{sec}.\arabic{punct}.  }}

\def\COUNTERS{\addtocounter{sec}{1}
              \setcounter{punct}{0}
          \setcounter{equation}{0}
          \setcounter{theorem}{0}
                  }

\newtheorem{theorem}{Theorem}[sec]
\newtheorem{proposition}[theorem]{Proposition}

\newtheorem{lemma}[theorem]{Lemma}

 \def\ov{\overline}
\def\wt{\widetilde}

 \newcommand{\rk}{\mathop {\mathrm {rk}}\nolimits}

 \newcommand{\indef}{\mathop {\mathrm {indef}}\nolimits}
 \newcommand{\dom}{\mathop {\mathrm {dom}}\nolimits}

     \newcommand{\Aut}{\mathop {\mathrm {Aut}}\nolimits}
     \newcommand{\Ams}{\mathop {\mathrm {Ams}}\nolimits}
     \newcommand{\Pol}{\mathop {\mathrm {Pol}}\nolimits}
     \newcommand{\Polh}{\mathop {\mathrm {Polh}}\nolimits}
 
\begin{document}

\def\OO{\mathrm{O}}
\def\GLO{\mathrm{GLO}}
\def\Coll{\mathrm{Coll}}
\def\kappa{\varkappa}
\def\Mat{\mathrm{Mat}}
\def\U{\mathrm U}
\def\GLL{\overline{\mathrm {GL}}}
\def\Spp{\overline{\mathrm {Sp}}}

\def\R{\mathbb{R}}
\def\C{\mathbb{C}}
\def\V{\mathbb{V}}

\def\la{\langle}
\def\ra{\rangle}

 \def\cA{\mathcal A}
\def\cB{\mathcal B}
\def\cC{\mathcal C}
\def\cD{\mathcal D}
\def\cE{\mathcal E}
\def\cF{\mathcal F}
\def\cG{\mathcal G}
\def\cH{\mathcal H}
\def\cJ{\mathcal J}
\def\cI{\mathcal I}
\def\cK{\mathcal K}
 \def\cL{\mathcal L}
\def\cM{\mathcal M}
\def\cN{\mathcal N}
 \def\cO{\mathcal O}
\def\cP{\mathcal P}
\def\cQ{\mathcal Q}
\def\cR{\mathcal R}
\def\cS{\mathcal S}
\def\cT{\mathcal T}
\def\cU{\mathcal U}
\def\cV{\mathcal V}
 \def\cW{\mathcal W}
\def\cX{\mathcal X}
 \def\cY{\mathcal Y}
 \def\cZ{\mathcal Z}
 
 \def\cGL{\mathcal{GL}}
 
 \def\EXP{\mathcal{EXP}}
 
 \def\frA{\mathfrak A}
 \def\frB{\mathfrak B}
\def\frC{\mathfrak C}
\def\frD{\mathfrak D}
\def\frE{\mathfrak E}
\def\frF{\mathfrak F}
\def\frG{\mathfrak G}
\def\frH{\mathfrak H}
\def\frI{\mathfrak I}
 \def\frJ{\mathfrak J}
 \def\frK{\mathfrak K}
 \def\frL{\mathfrak L}
\def\frM{\mathfrak M}
 \def\frN{\mathfrak N} \def\frO{\mathfrak O} \def\frP{\mathfrak P} \def\frQ{\mathfrak Q} \def\frR{\mathfrak R}
 \def\frS{\mathfrak S} \def\frT{\mathfrak T} \def\frU{\mathfrak U} \def\frV{\mathfrak V} \def\frW{\mathfrak W}
 \def\frX{\mathfrak X} \def\frY{\mathfrak Y} \def\frZ{\mathfrak Z} \def\fra{\mathfrak a} \def\frb{\mathfrak b}
 \def\frc{\mathfrak c} \def\frd{\mathfrak d} \def\fre{\mathfrak e} \def\frf{\mathfrak f} \def\frg{\mathfrak g}
 \def\frh{\mathfrak h} \def\fri{\mathfrak i} \def\frj{\mathfrak j} \def\frk{\mathfrak k} \def\frl{\mathfrak l}
 \def\frm{\mathfrak m} \def\frn{\mathfrak n} \def\fro{\mathfrak o} \def\frp{\mathfrak p} \def\frq{\mathfrak q}
 \def\frr{\mathfrak r} \def\frs{\mathfrak s} \def\frt{\mathfrak t} \def\fru{\mathfrak u} \def\frv{\mathfrak v}
 \def\frw{\mathfrak w} \def\frx{\mathfrak x} \def\fry{\mathfrak y} \def\frz{\mathfrak z} \def\frsp{\mathfrak{sp}}
 \def\bfa{\mathbf a} \def\bfb{\mathbf b} \def\bfc{\mathbf c} \def\bfd{\mathbf d} \def\bfe{\mathbf e} \def\bff{\mathbf f}
 \def\bfg{\mathbf g} \def\bfh{\mathbf h} \def\bfi{\mathbf i} \def\bfj{\mathbf j} \def\bfk{\mathbf k} \def\bfl{\mathbf l}
 \def\bfm{\mathbf m} \def\bfn{\mathbf n} \def\bfo{\mathbf o} \def\bfp{\mathbf p} \def\bfq{\mathbf q} \def\bfr{\mathbf r}
 \def\bfs{\mathbf s} \def\bft{\mathbf t} \def\bfu{\mathbf u} \def\bfv{\mathbf v} \def\bfw{\mathbf w} \def\bfx{\mathbf x}
 \def\bfy{\mathbf y} \def\bfz{\mathbf z} \def\bfA{\mathbf A} \def\bfB{\mathbf B} \def\bfC{\mathbf C} \def\bfD{\mathbf D}
 \def\bfE{\mathbf E} \def\bfF{\mathbf F} \def\bfG{\mathbf G} \def\bfH{\mathbf H} \def\bfI{\mathbf I} \def\bfJ{\mathbf J}
 \def\bfK{\mathbf K} \def\bfL{\mathbf L} \def\bfM{\mathbf M} \def\bfN{\mathbf N} \def\bfO{\mathbf O} \def\bfP{\mathbf P}
 \def\bfQ{\mathbf Q} \def\bfR{\mathbf R} \def\bfS{\mathbf S} \def\bfT{\mathbf T} \def\bfU{\mathbf U} \def\bfV{\mathbf V}
 \def\bfW{\mathbf W} \def\bfX{\mathbf X} \def\bfY{\mathbf Y} \def\bfZ{\mathbf Z} \def\bfw{\mathbf w}
 \def\R {{\mathbb R }} \def\C {{\mathbb C }} \def\Z{{\mathbb Z}} \def\H{{\mathbb H}}
  \def\K{{\mathbb K}}
   \def\k{{\Bbbk}}
 \def\N{{\mathbb N}} \def\Q{{\mathbb Q}} \def\A{{\mathbb A}} \def\T{\mathbb T} 
 \def\G{\mathbb G}
 \def\bbA{\mathbb A} \def\bbB{\mathbb B} \def\bbD{\mathbb D} \def\bbE{\mathbb E} \def\bbF{\mathbb F} \def\bbG{\mathbb G}
 \def\bbI{\mathbb I} \def\bbJ{\mathbb J} \def\bbL{\mathbb L} \def\bbM{\mathbb M} \def\bbN{\mathbb N} \def\bbO{\mathbb O}
 \def\bbP{\mathbb P} \def\bbQ{\mathbb Q} \def\bbS{\mathbb S} \def\bbT{\mathbb T} \def\bbU{\mathbb U} \def\bbV{\mathbb V}
 \def\bbW{\mathbb W} \def\bbX{\mathbb X} \def\bbY{\mathbb Y} \def\kappa{\varkappa} \def\epsilon{\varepsilon}
 \def\phi{\varphi} \def\le{\leqslant} \def\ge{\geqslant}

\def\P{\mathbf P}

\def\GL{\mathrm {GL}}
\def\bGL{\mathbf {GL}}
\def\GLB{\mathrm {GLB}}

\def\bGr{\mathbf {Gr}}
\def\Gr{\mathrm {Gr}}
\def\Sp{\mathrm {Sp}}
\def\bFl{\mathbf {Fl}}

\def\1{\mathbf {1}}
\def\0{\mathbf {0}}

\def\rra{\rightrightarrows}
\def\rt{\rightarrowtail}

 \newcommand{\Dim}{\mathop {\mathrm {Dim}}\nolimits}
  \newcommand{\codim}{\mathop {\mathrm {codim}}\nolimits}
   \newcommand{\im}{\mathop {\mathrm {im}}\nolimits}
\newcommand{\ind}{\mathop {\mathrm {ind}}\nolimits}
\newcommand{\graph}{\mathop {\mathrm {graph}}\nolimits}

\def\F{\bbF}

\def\lambdA{{\boldsymbol{\lambda}}}
\def\alphA{{\boldsymbol{\alpha}}}
\def\betA{{\boldsymbol{\beta}}}
\def\gammA{{\boldsymbol{\gamma}}}
\def\deltA{{\boldsymbol{\delta}}}
\def\mU{{\boldsymbol{\mu}}}
\def\nU{{\boldsymbol{\nu}}}
\def\epsiloN{{\boldsymbol{\varepsilon}}}
\def\phI{{\boldsymbol{\phi}}}
\def\psI{{\boldsymbol{\psi}}}

\def\pil{\overset{\leftarrow}{\pi}}
\def\pir{\overset{\rightarrow}{\pi}}

\def\sm{\smallskip}
\def\nw{\nwarrow}
\def\se{\searrow}

\begin{center}
\Large\bf

Polyhomomorphisms of locally compact groups

\medskip

\large \sc
Yury A. Neretin%
\footnote{Supported by the grant  FWF, P31591.}

\end{center}

\COUNTERS

{\small Let $G$ and $H$ be locally compact groups with fixed two-side-invariant Haar measures.
A polyhomomorphism  $G\rt H$ is a closed subgroup $R\subset G\times H$ with a fixed Haar measure, whose marginals
on $G$ and $H$ are dominated by the Haar measures on $G$ and $H$.
A polyhomomorphism can be regarded as a multi-valued map sending points to sets equipped with
'uniform' measures.
  For polyhomomorphsisms $G\rt H$, $H\rt K$ there is a well-defined product
 $G\rt K$. The set of polyhomomorphisms $G\rt H$ is a metrizable compact space
 with respect to the Chabauty--Bourbaki topology and the product is separately continuous.
 A polyhomomorphism $G\rt H$ determines a canonical operator
 $L^2(H)\to L^2(G)$, which is a partial isometry
 up to scalar factor.
 As an example, we  consider locally compact infinite-dimensional linear spaces over
 finite fields and examine  closures of groups of linear operators
 in semigroups of polyendomorphisms.  
}

\bigskip

\section{Polyhomomorphisms}

{\it We consider only second-countable  locally compact groups} $G$, i. e., locally group having
 countable base of open subsets. Such groups admit left invariant metrics (see, e.g. \cite{HR}, Theorem 8.3),
so they are separable as metric spaces. They are complete topological groups  in the sense of Weil--Bourbaki
(\cite{Bour1}, Corollary III.3.1). Such a group admits a unique up to a scalar factor
left-invariant measure $\gamma(g)$ (the Haar measure), see  \cite{Bour2}, Theorem VII.1.1, \cite{HR}, Theorem 8.3, \cite{Zhe}, Theorem 26.4, \cite{Wei}, \cite{Ner-Haar}, Chapter A.

Since $G$ has a structure of complete metric space, the Borel structure on $G$ is standard 
(see, e.~g., \cite{Kech}, Sect 12), as a Borel space $G$ is isomorphic to the line $\R$, 
or a  countable set, or a finite set. As a space with measure
the group $G$ is a Lebesgue space (see, e.~g., \cite{Bog}, Chapter 10), as a space with 
measure  $G$ can be isomorphic
to an interval, the line $\R$, a  countable or finite set .

Since we have a measure, we also have standard  spaces of measurable functions on $G$ as $L^2(G)$. 
By $C_c(G)$ we denote the space of continuous functions on $G$ with compact support.

A locally compact group is {\it unimodular} if the Haar measure is two-side invariant
(see \cite{Bour2}, Subsect. VII.1.3-4).

\sm

Let $K\supset L$ be groups. Denote by $[K:L]$ the {\it index}, i.e., the number
of elements in $K/L$.

\sm 

For a set $X$ and subset $A\subset X$ we denote by $I_A(x)$ the {\it indicator function} of $A$, i.e.,
$I_A(x)=1$ if $x\in A$ and $0$ otherwise.

\sm

{\bf \punct Multiplicative relations.} Let $X$, $Y$ be sets. 
A {\it relation $X\rra Y$} is a subset $R\subset X\times Y$. For two relations
$R:X\rra Y$, $S:Y\rra Z$ we define their {\it product} $SR:X\rra Z$ as the set of all
$(x,y)\in X\times Z$, for which there exists $y\in Y$ such that $(x,y)\in R$, $(y,z)\in S$.
Clearly, this product is associative.

For a relation $R:X\rra Y$ we define: 

\sm

--- the {\it image} $\im R$ is the projection of $R$ to $Y$;

\sm

--- the {\it domain} $\dom R$ is the projection of $R$ to $X$.

\sm

Define the {\it pseudoinverse} relation $R^\square:Y\rra X$ as the same subset 
 $R\subset X\times Y$ considered as a relation
from $Y$ to $X$. Obviously,
$$
(TR)^\square=R^\square T^\square.
$$

For a subset $A\subset X$ we define {\it its image} $RA$ as the set of all $b\in Y$ such that there exists
$a\in A$ satisfying $(a,b)\in R$. 

\sm 

{\sc Remark.} If $f:X\to Y$ is a map, then its graph $\Gamma(f)\subset X\times Y$
is a relation, $\dom \Gamma(f)=X$  and the projection map $\Gamma(f)\to X$ is injective.
\hfill $\boxtimes$

\sm

A {\it partial bijection} $X\to Y$ is a bijective map of a subset $A\subset X$ to a subset $B\subset Y$.
A relation $R:X\rra Y$  is a partial bijection if the projection maps
from $R$ to $X$ and $Y$ are injective.

\sm

Let $G$, $H$ be groups. A {\it multiplicative relation} $R:G\rra H$ is a subgroup
in $G\times H$. Clearly, a product of multiplicative relations is a multiplicative relation.
For a multiplicative relation $R:G\rra H$ we define

\sm

--- the {\it kernel} as the intersection of $R$ with $G\subset G\times H$;

\sm

--- the {\it indefinity} as the intersection of $R$ with $H\subset G\times H$.

\sm

The following statement is obvious.

\begin{lemma}
	\label{l:lemma}
{\rm a)} The kernel $\ker R$ is a normal subgroup in $\dom R$ and $\indef R$ is a normal subgroup in $\im R$.

\sm 

{\rm b)} A multiplicative relation $R$ determines a canonical isomorphism
$$
\iota (R):\dom R/ \ker R\to \im R/\indef R.
$$

{\rm c)} The subgroups $\ker R$, $\indef R$, $\ker R\times \indef R$ are normal in $R$.
\end{lemma}

\sm

We define a {\it partial isomorphism} $G\to H$ as a partial bijection between two subgroups $A\subset G$,
$B\subset H$ sending products to products.

\sm

If groups $G$, $H$ are additive, then it is reasonable to say '{\it additive relation}'.
If they are linear spaces and $R$ is a subspace, we say '{\it linear relation}'.
Linear relations and additive relation are usual mathematical objects (see, e.g., \cite{Mac}, \cite{Ner-gauss}), multiplicative relations are known but
appear not too often, see, e.g., \cite{Sch}, Sect. 1.2.

\sm 

{\bf \punct The category of polyhomomorphisms.}
Let $X$ be a space with measure $\xi$, let $Y$ be a set
and $f$ be a map $X\to Y$. Recall that the {\it image $\upsilon$ of the measure} $\xi$ under
the map $f$ is the measure on $Y$ defined by the condition
$\upsilon(B):=\xi\bigl(f^{-1}(B)\bigr)$.

Denote by $G^\circ$ (resp. $H^\circ$) a unimodular group $G$ (resp. $H$) with a fixed Haar measure $\gamma(g)$ (resp. $\eta(h)$).
Denote by $\pil$ the natural projection map $G\times H\to G$, by $\pir$ the projection map
$G\times H\to H$
We say that a {\it polyhomomorphism} $R^\circ:G^\circ \rt H^\circ$
is an object of one of the following types:

\sm

1. a closed
subgroup $R\subset G\times H$ with a fixed Haar measure 
$\rho(r)$ such that the image of $\rho$ under  $\pil$ (respectively, $\pir$) is
dominated by $\gamma(g)$ (resp. by $\eta(h)$);

\sm

2. the zero measure $0=0_{G,H}$ on $G\times H$.

\sm

Denote by $\Polh(G^\circ,H^\circ)$ the set of all polyhomomorphisms $G^\circ\rt H^\circ$. 
Elements of this set automatically satisfy the following properties
(so they can be included to the definition of polyhomomorphisms).

\begin{proposition}
\label{pr:add-def}
Let $R^\circ\in \Polh(G^\circ,H^\circ)$ and $R$ the underlying multiplicative relation. Then

\sm

 {\rm a)} The subgroups $\ker R\subset G$, $\indef R\subset H$
 are compact. 
 
 \sm
 
 {\rm b)} The subgroups $\dom R\subset G$, $\im R\subset H$
 are open
 
 \sm
 
{\rm c)} The group $R$ is unimodular.

\sm 

{\rm d)}
The image of the measure
$\rho(r)$ under the projection $\pil:R\to G$ 
is the measure on $\dom G$ having the form $\alpha\, \gamma(g)$,
where
$$\alpha=\alpha(R^\circ)$$
is a constant such that
  $0<\alpha\le 1$. Similarly,
the image of $\rho(r)$ under  $\pir $ is the measure on
 $\im R$ of the form  $\beta\, \eta(h)$, where
$$\beta=\beta(R^\circ)$$
is a constant satisfying the condition
$0<\beta\le 1$.
\end{proposition}

The proof is contained in Subsect. \ref{ss:prop1}.

\sm

Next, let $R^\circ\in (G^\circ, H^\circ)$, $T^\circ\in \Polh(H^\circ, K^\circ)$
be two nonzero polyhomomorphisms. We define their {\it product}
$S^\circ= T^\circ R^\circ\in \Polh(G^\circ, K^\circ)$ as follows:

\sm

1) a multiplicative relation $S$ is $S:=TR$;

\sm

2) we normalize the Haar measure on $S$ in the terms of its images under the projections to $G$ and $K$:
\begin{align}
\alpha(S^\circ)&=\frac{\alpha(R^\circ)\, \alpha(T^\circ)}
{ \bigl[\indef R:(\indef R\cap \dom T)\bigr]}
;
\label{eq:product1}
\\
\beta(S^\circ)&=\frac{\beta(R^\circ)\,\beta (T^\circ)}
{\bigl[\ker T:(\ker T\cap \im R)\bigr]}.
\label{eq:product2}
\end{align}

\sm

A product of a zero polyhomomorphism and any  polyhomomorphism is zero,
$$
0_{H,K} R^\circ=0_{G,K}, \qquad T^\circ\, 0_{G,H}=0_{G,K}, \qquad
0_{H,K}0_{G,H} =0_{G,K}.
$$

\begin{lemma}
\label{l:ass}
	This product is well defined and associative, i.~e., for any $G^\circ$,
	$H^\circ$, $K^\circ$, $L^\circ$ and any
	$$
	R^\circ\in \Polh(G^\circ,H^\circ),\quad T^\circ\in \Polh(H^\circ,K^\circ),
	\quad S^\circ\in \Polh(K^\circ,L^\circ)
	$$
	we have
	$$
	(S^\circ T^\circ)R^\circ=S^\circ(T^\circ R^\circ).
	$$
\end{lemma}

\sm

{\sc Remark.}
In fact the definition of a product becomes a theorem if we consider polyhomomorphisms as 
special cases of polymorphisms, see below Subsect. \ref{ss:polymorphisms}--\ref{ss:poly-poly-homo}.
\hfill $\boxtimes$

\sm

Thus we get a {\it category of polyhomomorphisms}. Objects $G^\circ$ of this category
are unimodular locally compact groups equipped with  fixed Haar measures.
The set of morphisms from $G^\circ$ to $H^\circ$ is
$\Polh(G^\circ,H^\circ)$.

\sm 

For $R^\circ\in\Polh(G^\circ,H^\circ)$ denote by $(R^\circ)^\square$ the same subgroup 
in $G\times H$ with the same Haar measure considered as a subgroup in $H\times G$.
Obviously, $(T^\circ R^\circ)^\square=(R^\circ)^\square (T^\circ)^\square$.

\sm

{\bf \punct Comments to the definition.\label{ss:finite}}
Here we briefly discuss, what means this definition
for some natural classes of locally compact groups.

\sm

{\sc Finite groups.} Let groups  $G$, $H$ be finite.
Normalize measures assuming that each element has measure 1.
A subgroup 
$R\subset G\times H$ can be arbitrary, we must equip it with a uniform measure.
Assuming that a measure of each element of the group  $R$ is $\frr$,
we get the following inequalities for  $\frr$:
$$
\alpha(R^\circ) =\frr\cdot \# \indef R\le1,\qquad \beta(R^\circ)= \frr\cdot \#\ker R\le 1,
$$
the symbol $\# X$ denotes the number of elements of a set  $X$.

The category of polyhomomorphisms is a kind of a 'central extension'
of the category of multiplicative linear relations%
\footnote{See a formal definition in   \cite[Subsect.~I.8.6]{Ner-cl}.
	Notice that categories of linear relations arising in representation theory
	appear together with their central extensions, see
	\cite{Ner-cl}, \cite{Ner-book}.}. Let explain an origin of this extension
with an example of  finite groups. Denote by  $\ell^2(G)$ the space of functions on   $G$
equipped with the  $\ell^2$-inner product. Notice that any homomorphism 
$\rho:G\to H$ determines an operator 
$$\Pi_*(\rho): \ell^2(H)\to \ell^2(G)$$
by the formula 
$$
\Pi_*(\rho) f(g):=f(\rho(g))
.$$ 
If $\sigma$ is a homomorphism  $H\to K$, then
$$
\Pi_*( \rho\circ \sigma)=\Pi_*(\rho) \Pi_*(\sigma).
$$ 
For a multiplicative relation
$R:G\rra H$ we can define the following operator  $\ell^2(H)\to \ell^2(G)$:
\begin{equation}
\Pi_*(R)f(g):=\sum_{h\in H:\, (g,h)\in R} f(h).
\label{eq:sum}
\end{equation}
It easy to show that for multiplicative relations
$R:G\rra H$, $T:H\rra K$ the following identity holds 
\begin{multline*}
\Pi_*(R)\,\Pi_*(T)=\#(\ker T\cap \indef R)\cdot  \Pi_*(TR)
=\\=
\frac{\# \indef T\cdot \#(\indef R\cap \dom T)}{\#\indef TR} \Pi_*(TR)=
\frac{\#\ker R\cdot \#(\ker T\cap \im R)}{\#\ker TR}
\Pi_*(TR).
\end{multline*}

As above, define a uniform measure
on $R$ and modify operators assuming
$$
\Pi(R^\circ):= \frr \cdot \Pi_*(R).
$$
It is easy to see that
$$
\Pi(R)\,\Pi(T)=  \Pi(T^\circ R^\circ).
$$
An extension of this construction to arbitrary polyhomomorphisms
is discussed below in Subsect.
\ref{ss:homothety}.

\sm

{\sc Discrete groups.} If $G$, $H$ are discrete   (countable),
then we have an additional condition for 
$R$: the groups  $\ker R$ and $\indef R$ must be finite.
Notice also, that this condition is 
necessary and sufficient for the boundedness of the operator 
(\ref{eq:sum}).

\sm 

{\sc Lie groups.} For connected Lie  groups the notion of polyhomomorphism gives little new
in comparation with homomorphisms. By definition, connected groups do not have proper 
open subgroups, on other hand there are few compact normal subgroups in Lie groups.
Therefore it remains little possibilities 
to satisfy conditions of  Proposition
\ref{pr:add-def}.a-b and
Lemma \ref{l:lemma}.a. 

\sm 

For a reader familiar with theory of Lie groups,
we present a typical 
{\it example of a polyhomomorphism between semisimple groups}.
Consider the group  $G:=\mathrm{SL}(2,\R)$ (the group of real matrices of order
2 with determinant  1, it is homotopically equivalent
to the circle, so its fundamental group is
$\Z$). Denote by   $G^{\sim n}$ its $n$-sheeted covering. The group $G^{\sim 6}$
embeds to  $G^{\sim 2}\times G^{\sim 3}$ and can be regarded as a multiplicative
relation, it remains to normalize a Haar measure on
$G^{\sim 6}$ in some way.

\sm

{\it The case of tori} gives more possibilities.
Recall that a torus is a quotient 
$\R^n/ \Z^n$, we equip it with a probabilistic Haar measure.
We need closed subgroups in $\R^n/ \Z^n\times \R^m/ \Z^m$.
They can be easily described. Namely let $L$ be a linear subspace
in $\R^n\times\R^m$ determined by equations with integer coefficients,
$$
\Bigr\{\sum_{i=1}^n p^i_\alpha x_i+\sum_{j=1}^m q^j_\alpha y_j=0,
\qquad \text{where $\alpha=1$, \dots, $k$,} 
$$
or, in a matrix form
$Px+Qy=0$. The image of $L$ under the map
$\R^n\times\R^m\to \R^n/ \Z^n\times \R^m/ \Z^m$ is a closed subgroup
in the product of tori, and all closed subgroups have such a form.  

Next, consider a finitely generated subgroup
$\Gamma$ in a linear space $\Q^k$ over rationals.
Take the set   $L_\Gamma$ of vectors in $\R^{n+m}$,
satisfying the condition
$Px+Qy\in \Gamma$ and consider its projection to  
$ \R^n/ \Z^n\times \R^m/ \Z^m$. This construction gives us 
all closed subgroups in a product of tori. 

To get a polyhomomorphism, we need the surjectivity of the projections
of
$L$ to $\R^n$ and $\R^m$. Notice that in this case
the denominators in formulas  \eqref{eq:product1}-\eqref{eq:product2}
equal  1. Therefore nothing prevents us without loss of substance
set $\alpha=\beta=1$ for all polyhomomorphisms.

\sm

{\sc Totally disconnected non-discrete groups.}
There are lot of totally disconnected locally compact groups arising
in various branches of mathematics (for instance,  $p$-adic and adelic groups,
Galois groups equipped with the Krull topology, groups of automorphisms and spheromorphisms
of trees, some infinite-dimensional groups over finite fields.

Notice that each group 
$G$ of this type has 'many' {\it polyendomorphisms} (i.e., polyhomomorphisms to itself).
We present mass examples:

\sm

1) Let $\Omega\subset G$ be an open subgroup.
We take the subgroup  $R\subset G\times G$ consisting of points
of the form,  $(g,g)$, where  $g$ ranges in  $\Omega$.
Evidently,
$R\simeq \Omega$, and this determines a Haar measure on   $\Omega$.

\sm

2) Take an open subgroup  $\Omega\subset G$ and a normal compact open subgroup%
\footnote{A totally disconnected locally compact group has a fundamental system
	of neighborhoods of the unit consisting of open compact subgroup. If the group
	is compact, then these subgroups can be chosen normal,
	see~\cite[Theorems 7.5-7.7]{HR}. This implies that there is a countable
	set of open subgroups in  $\Omega$; for a fixed open compact
	subgroup  $\Omega$ we have a countable set of open normal subgroups.} $L\subset \Omega$.
We take $R\subset G\times G$ consisting of elements
of the form 
$$(l g,g),\qquad\text{where $g\in \Omega$, $l\in L$.}
$$
It remains to equip 
$R$ with a Haar measure  (and multiply it by a sufficiently
small factor if this is necessary).

\sm 

3) Any automorphism $\theta$ of a group  $G$ (for instance, an interior
automorphism) generates a polyhomomorphism  $G^\circ\rt  G^\circ$.
As $R$ we take a graph of the map
$\theta$ and equip it with an appropriate Haar measure%
\footnote{Generally speaking, the measures $\gamma(g)$ and $\gamma(\theta(g))$
	differs by a scalar factor $c$, see. e.g., \cite{Bour2}, Subsect. VII.1.4. If $c\ne 1$, then the coefficients $\alpha(R^\circ)$
	and $\beta(R^\circ)$ are different.}. 

\sm

Further, we can multiply polyhomomorphisms of types
1)-3).

\sm

{\bf\punct The convergence of polyhomomorphisms.%
\label{ss:ph-convergence}}
We define a convergence in $\Polh(G^\circ,H^\circ)$
as the weak convergence of measures (see, e.~g., \cite{Bog}, Sect.8.1) on $G\times H$.
Namely, let $R_j^\circ=(R_j,\rho_j)$, $R^\circ=(R,\rho)$
be polyhomomorphisms $G\rt H$.
A sequence $R^\circ_j:G^\circ\rt H^\circ$ {\it converges} to 
$R^\circ:G^\circ\rt H^\circ$  if for any functions $\phi\in C_c(G)$ and $\psi \in C_c(H)$
we have the converges of integrals 
$$
\int_{R_j} \phi\bigl(\pil(r)\bigr)\,\psi\bigl(\pir(r)\bigr)\,d \rho_j(r)
\to\int_{R} \phi\bigl(\pil(r)\bigr)\,\psi\bigl(\pir(r)\bigr)\,d \rho(r).
$$
A sequence $R^\circ_j$ converges to 
$0_{G,H}$ if for any $\phi$, $\psi$ we have the convergence
$$
\int_{R_j} \phi\bigl(\pil(r)\bigr)\,\psi\bigl(\pir(r)\bigr)\,d \rho_j(r)\to 0.
$$

{\sc Remark.}
We can define the convergence in the following equivalent way:
$R_j^\circ$ converges to $R^\circ$ if
 for any $\theta\in C_c(G\times H)$
we have the convergence
$$
\int_{R_j} \theta(r)\,d \rho_j(r)\to \int_{R} \theta(r)\,d\rho(r).
$$
Similarly, $R_j^\circ$ converges to $0_{G,H}$ if such sequences
of integrals converge to 0. \hfill $\boxtimes$

\begin{proposition}
\label{pr:2}
 {\rm a)} This convergence is metrizable and sets $\Polh(G^\circ,H^\circ)$ are 
  compact.
 
 \sm
 
 {\rm b)} The product of polyhomomorphisms is separately continuous.
\end{proposition}

This convergence  is a rephrasing
the Chabauty--Bourbaki convergence of subgroups in locally compact groups,
 see Bourbaki  \cite{Bour2}, Sect VIII.5,
Bourbaki normalizes Haar measures on each subgroup, we allow to vary scalar factors.
The compactness is  Theorem VIII.5.1 of Bourbaki.

\sm 

{\sc Remark.} Convergence $R^\circ_j\to R^\circ$
implies a convergence of subsets $R_j\to R$.
There are many non-equivalent definitions
of convergences on  sets of closed subsets of  topological or metric spaces,
see, e.~g., \cite{Mich}.
Our space $G\times H$ is locally compact and reasonable topologies coincide.  
For instance (see Bourbaki \cite{Bour1}, Subsect. VIII.5.6) we can take a left invariant metric on $G\times H$ compatible with the topology
and say that a sequence $R_j$ of closed subgroups converges to $R$ if for each $\epsilon>0$
for any compact set $K\subset G\times H$ for sufficiently large $j$ the set $K\cap R$ is contained in the 
$\epsilon$-neighborhood of $R_j$ and $K\cap R_j$ is contained in the $\epsilon$-neighborhood of
of $R$. See \cite{Bir} on a way to define a metric on this space.
 \hfill $\boxtimes$

\sm

{\sc Remark.} A convergence $R^\circ_j\to R^\circ$
does not implies convergences $\ker R_j\to \ker R_j$, $\dom R_j\to \dom R$,
$\alpha(R_j^\circ) \to \alpha(R_j^\circ)$, etc. However, we have some semicontinuities.
If $\ker R_j$ contain some subgroup $L\subset G$ starting some $j$, then
$\ker R$ contains $L$. If $\dom R_j$ are contained in a certain subgroup
$M\subset G$ starting some $j$, then $\dom R$ is contained in the same subgroup. 
If $\alpha(R^\circ_j)\le s$ starting some $j$, then $\alpha(R^\circ)\le s$.
\hfill $\boxtimes$

\sm

{\bf \punct Polymorphisms. Preliminaries.%
\label{ss:polymorphisms}}
See \cite{Ner-bist}, \cite{Ner-book}, Sect VIII.4.

\sm 

{\sc A. Category of polymorphisms.} Recall that a space $X$ with a finite or $\sigma$-finite 
measure $\xi$
is a
{\it Lebesgue measure space}   if it is equivalent to a union
of a finite or infinite interval of the line $\R$ and of a finite (may be, empty) or countable
collection of points having positive measures.

\sm

{\sc Remark.}
A locally compact group $G$ equipped with
the Haar measure as a measure space is equivalent
to

\sm 

--- a collection of points having equal positive measures if a group is discrete;

\sm

--- a finite interval $(a,b)\subset \R$ if a group is compact and infinite;

\sm

--- $\R$ otherwise.
\hfill$\boxtimes$

\sm

Let $(X,\xi)$ and $(Y,\upsilon)$
be Lebesgue  spaces.
A {\it polymorphism} $\mu:(X,\xi)\rt(Y,\upsilon)$ is a measure $\alpha$
on $X\times Y$ such that projection of $\mu$ to $X$ is dominated by $\xi$
(i.e., for any subset $A\subset X$ of finite measure we have $\xi(A)\ge \mu(A\times Y)$)
and the projection of $\mu$ to $Y$ is dominated by $\upsilon$.
We admit zero measures.
Denote the set of all polymorphisms $X\rt Y$ by $\Pol(X,Y)$.

We regard a polymorphism as a 'multivalued maps' $X\to Y$. Namely, for any
polymorphism $\mu:X\rt Y$ 
there is a canonical map (defined a.~s.) sending  points $x\in X$ to conditional measures
(see, e.g., \cite{Bog}, Sect. 10.4) $\mu_x(y)$ on $Y$ such that
for any subsets $A\subset X$, $B\subset Y$ of finite measure we have
$$
\mu(A\times B)=\int_A \mu_x(B)\,d\xi(x).
$$
Notice that  $\mu_x(Y)\le 1$ for almost all $x\in X$, also
we have
$$
\int_X \mu_x(B)\,d\xi(x)\le \mu(B).
$$

 Let $\mu:X\rt Y$, $\nu:Y\rt Z$ be polymorphisms.
 We define their {\it product} $\kappa=\nu\mu:X\rt Z$
 in the terms of conditional measures:
 $$
 \kappa_x=\int_Y \nu_y\,d\mu_x(y).
 $$
Thus we get a category whose objects are Lebesgue measure spaces and morphisms
are polymorphisms. 

For a polymorphism $\mu:X\rt Y$
we define the {\it adjoint
polymorphism} $\mu^\square:Y\to X$ that is the same measure considered as a measure on 
$Y\times X$.

\sm

{\sc B. Linear operators determined by polymorphisms.}
For a polymorphism $\mu:X\rt Y$ we consider the sesquilinear
form 
$$S_\mu:L^2(X,\xi)\times L^2(Y,\upsilon)\to \C$$
 defined by
$$
S_\mu(\phi,\psi)=
\int_{X\times Y} \phi(x)\,\ov {\psi(y)}\, d\mu(x,y).
$$
Applying the Cauchy--Bunyakovsky inequality and the definition of polymorphisms
we get 
\begin{equation}
|S_\mu(\phi,\psi)|\le \|\phi\|_{L^2(X,\xi)}\cdot\|\psi\|_{L^2(Y,\upsilon)}.
\label{eq:Bun}
\end{equation}
Therefore there exists a bounded  operator 
$$\Pi(\mu): L^2(Y,\upsilon)\to L^2(X,\xi)$$
 such that
for all $\phi\in L^2(X,\xi) $, $\psi\in L^2(Y,\upsilon)$ we have
$$
S_\mu(\phi,\psi)=\la \phi, \Pi(\mu) \psi\ra_{L^2(X,\xi)}.
$$
By (\ref{eq:Bun}), operators $\Pi(\mu)$ are {\it contractions}, i.~e.,
$$
\|\Pi(\mu)\|\le 1.
$$
The explicit expression for this operators is
\begin{equation}
\Pi(\mu)\psi(x)=\int_Y \psi(y)\,d\mu_x(y),
\label{eq:Pi}
\end{equation}
where $\mu_x$ are the conditional measures defined above.

\sm 

{\sc Remark.} 
The last expression shows that $\Pi(\mu)$ sends nonnegative functions to nonnegative functions.
Conversely, say that  an operator $T:L^2(Y,\xi)\to  L^2(X,\xi)$ is a {\it sub-Markov operator}
if it satisfies this property
and $\|T\|\le 1$. It is easy to show that any
sub-Markov operator has the form $T=\Pi(\mu)$ for some polymorphism $\mu$.
The measure $\mu$ is determined by the condition
$$
\mu(A\times B):=\la I_A, T I_B\ra_{L^2(X,\xi)},
$$
where $A\subset X$, $B\subset Y$ are sets of finite measure.
\hfill $\boxtimes$

\sm

{\sc Remark.} We also can describe  the operator $\Pi(\mu)$
 in the following way. Consider  a bounded nonnegative function 
$\psi$ on $Y$ and the measure 
\begin{equation}
\psi(y)\,\mu(x,y)
\label{eq:psimu}
\end{equation}
on $X\times Y$. Taking its projection
to $X$ we get a measure, say $\Phi$, on $X$.
For each measurable subset $A\subset X$ we have 
$\Phi(A)=S_\mu(I_A,\psi)$, where $I_A$ is the indicator function.
Clearly for a set $C\subset X$ of zero measure we have $\Phi(C)=0$.
Therefore, the measure $\Phi$ is absolutely continuous with respect
to $\xi$, and we can define
$\Pi(\mu)\psi$ as the Radon--Nikodym derivative $d\Phi/d\xi$,
for  formulas for Radon--Nikodym derivatives, see, e.g., \cite{ShG}, Sect. 10.
\hfill $\boxtimes$

\sm

Formula (\ref{eq:Pi}) easily implies that
$$
\Pi(\nu \mu)=\Pi(\nu)\Pi(\mu).
$$
So we get a functor from the category of polymorphisms to the category
of Hilbert spaces and bounded operators.

\sm

{\sc C. Topology on sets $\Pol(X,Y)$.}
Next, let $\mu_j$, $\mu$ be polymorphisms $X\rt Y$.
We say that $\mu_j$ converges to $\mu$ if for any subsets
$A\subset X$, $B\subset Y$ of finite measure the sequence
$\mu_j(A\times B)$ converges to $\mu(A\times B)$.
It is easy to show that this convergence is equivalent to the
weak operator convergence of the corresponding sub-Markov operators $\Pi(\mu_j)\to \Pi(\mu)$, i.~e.,
$$
\la \phi,\Pi(\mu_j)\psi\ra_{L^2(X,\xi)}\to \la \phi,\Pi(\mu)\psi\ra_{L^2(X,\xi)}
\quad \text {for all $\phi\in L^2(X,\xi)$, $\psi \in L^2(Y,\upsilon)$.}
$$

Let $H$, $K$ be separable Hilbert spaces.
Denote by
 $\cC(H,K)$ the set of all contractions $H\to K$,
equip it with the weak operator topology.  
This set is compact metrizable
and the multiplications  
$$\cC(H,K)\times \cC(K,L)\to \cC(H,L)$$
 are separately continuous.
This easily implies that {\it sets $\Pol(X,Y)$ are compact metrizable
and the product is separately continuous}.

\sm

{\sc D. Measure preserving transformations and polymorphisms.}
 Let $(X,\xi)$ be a space with a $\sigma$-finite non-atomic measure
(i.e.,  let $X$ be equivalent to $\R$). Denote by $\Ams(X)$ the group of measure preserving transformations
of $X$. Let $g\in\Ams(X)$.  Consider the map $X\to X\times X$ given by
$x\mapsto (x,g(x))$, consider the image $\kappa_g$ of the measure $\xi$ under this map. It is clear 
that $\kappa_g$ is a polymorphism $X\rt X$ and products of measure preserving transformations
correspond to products polymorphisms. Also it is easy to show that
{\it the group $\Ams(X)$ is dense in the semigroup $\Pol(X,X)$}.

So the notion of a polymorphisms extends the notion of measure preserving transformation (apparently, this idea arises to E.~Hopf~\cite{Hopf}).

 \sm
 
 {\sc E. References on polymorphisms.} Polymorphisms and
 various variants of 'Markov operators' are the standard objects of ergodic theory, see, e.g.,
 E.~Hopf~\cite{Hopf}, J.~Neveu~\cite{Nev},
 A.\,V.~Vershik \cite{Ver}, U.~Krengel \cite{Kre}.
 There are several natural groups of transformation of measure spaces 
 (measure preserving transformations of spaces with finite measures or with $\sigma$-finite measures,
 groups of regular transformations, etc.). 
  For this reasons there are several kinds of 'polymorphisims', see \cite{Ner-bist},
 \cite{Ner-book}, Sect. VIII.4 and Chapter X. The version discussed above
 corresponds to the group of measure preserving transformations of a space with infinite
 continuous measure, apparently it
 appeared in \cite{Ner-bist}.
 
 K.~Schmidt and  A.\,V.~Vershik \cite{SV} considered polyhomomorphisms
 ('algebraic polymorphisms') of compact
 groups $K$ under stronger conditions.
 In our terminology they consider polyhomomorphisms $R^\circ:K^\circ\to K^\circ$
 such that $\dom R=K$, $\im R=K$ (also $\alpha(R^\circ)=\beta(R^\circ)=1$, but the last condition
 in this case
 is not essential). 
 In particular, this includes the case of tori discussed above in Subsect. \ref{ss:finite}.

 \sm

 {\bf\punct Polymorphisms and polyhomomorphisms.%
 \label{ss:poly-poly-homo}}
 So  {\it any polyhomomorphism is a polymorphism}.
 
 \begin{theorem}
 \label{th:product}
  The product of polyhomomorphism defined above corresponds to the product of polymorphisms.
 \end{theorem}

The proof occupies Subsect. \ref{ss:2support}--\ref{ss:2products-poly}.
 
 \sm
 
 This immediately implies the associativity of the product of polyhomomorphisms (Lemma \ref{l:ass}) and the separate continuity (Proposition \ref{pr:2}.b).

 \sm
 
 {\bf\punct Linear operators determined by polyhomorphisms.%
 \label{ss:homothety}}
  Let $G$ be a locally compact group, $\Phi$ an open subgroup, and $\Delta$ a  compact normal
 subgroup in
 $\Phi$. Denote by $L^2(\Phi)^\Delta\subset L^2(G)$ 
 the subspace consisting of functions that are supported by $\Phi$ and invariant with respect to
 $\Delta$. Denote by $P_{\Phi|\Delta}^G$ the operator of orthogonal projection
 to the subspace $L^2(\Phi)^\Delta\subset L^2(G)$.

 \begin{proposition}
 	\label{pr:homothety}
 	Let $R^\circ\in \Polh(G^\circ, H^\circ)$.
 	Then the operator
 	$$
\wt \Pi(R^\circ):= 	\bigl(\alpha(R^\circ)\,\beta(R^\circ)\bigr)^{-1/2}\, \Pi(R^\circ):\, L^2(H)\to L^2(G)
 	$$
 	is a partial isometry%
 	 	\footnote{Consider Hilbert spaces $V$, $W$.
 	 		A bounded linear operator   $A:V\to W$ 
 	 		is called a  {\it partial isometry},
 	 		if there exists a subspace $L\subset V$ ({\it the initial subspace}) and $M\subset W$ ({\it the final subspace}) such that the restriction of
 $A$ to $V$ is a unitary operator  $L\to M$
 and restriction of  $A$ to the orthogonal complement 
   $L^\bot\subset V$ is zero.}. The initial subspace
 	of $\wt \Pi(R^\circ)$ is $L^2(\im R)^{\indef R}$,  the final subspace is
 	$L^2(\dom R)^{\ker R}$.
 \end{proposition}

Explicit description of operators $\Pi(R^\circ)$ is contained in
Subsect. \ref{ss:2operators}.

\sm 

So operators $\Pi(R^\circ)$ are 'partial homotheties'.
By Theorem \ref{th:product}, a product $\Pi(T^\circ)\Pi(R^\circ)$ of two 'partial homotheties'
 is a 'partial homothety' again, but the new scaling coefficient is not a product
 of scaling coefficients. Formulate a geometric statement related to this phenomenon.

 Let $L$, $M$ be two closed subspaces in a Hilbert space $H$.
 Consider the operators $P_L$ and $P_M$ of orthogonal projections  
 $H\to L$, $H\to M$. Consider the self-adjoint operators
 $$
 P_L P_M\Bigr|_{L}=P_L P_M P_L\Bigr|_{L}: L\to L,\qquad P_MP_L\Bigr|_{M}=P_MP_L P_M\Bigr|_{M}:M\to M.
 $$
It is easy to see that their spectral types  coincide upto 
multiplicities of zeros, and these spectral types are invariants
of a pair of subspaces under unitary transformations (this is an analog of
angles in elementary geometry, see, e.g., \cite{Ner-cl}, Sect. 2.5).

 \begin{proposition}
 	\label{pr:angles}
 	Let $G$ be a unimodular locally compact group, $\Phi$, $\Psi$ 
open subgroups, $\Delta$ be a compact normal subgroup in $\Phi$, $\Gamma$ a compact
normal subgroup in $\Psi$. Let the pairs $(\Phi,\Delta)$ and $(\Psi,\Gamma)$
be different.  Then the spectrum of the operator
$$P_{\Phi|\Delta}^G P_{\Psi|\Gamma}^G\Bigr|_{L^2(\Phi)^\Delta}: L^2(\Phi)^\Delta \to  L^2(\Phi)^\Delta$$
consists of two points%
\footnote{This property is similar to isoclinicity.
	A pair of subspaces 
 $L$, $M$ of a finite dimensional Euclidean space is called {\it isoclinic},
 if the operator
	$P_LP_M:L\to M$ is scalar (or, equivalently,  all Jordan angles
	 between $L$, $M$ are equal), see, e.g. 
	\cite{Wolf}.}, namely, $0$ and 
$$
\sigma:=\bigl([\Delta:(\Delta\cap \Psi)]  \cdot [\Gamma:(\Gamma\cap \Phi)]\bigr)^{-1}
.$$
 Therefore the operator
$$
\sigma^{-1/2}P_{\Psi|\Gamma}^G\Bigr|_{L^2(\Phi)^\Delta}: L^2(\Phi)^\Delta\to L^2(\Psi)^\Gamma
$$
is a partial isometry.
 \end{proposition}
 
The proof is contained in \ref{ss:angles-proof}.

 \sm 
 
 {\bf \punct Rational polyhomomorphisms.}
 Let $K_1$, $K_2\subset G$ be open compact subgroups in a locally compact group $G$.
 Then $K_1\cap K_2$ also is an open compact subgroup. The homogeneous
 space $K_1/(K_1\cap K_2)$ is discrete and compact, therefore it is finite. 
 Therefore the ratio of measures of $K_1$ and $K_2$ is rational.
 
 Now we define a subcategory $\Polh_\Q$  of the category of polymorphisms. 
 Let us consider only unimodular locally compact groups
 that have open compact subgroups
 and equip them with Haar measures such that measures of open compact subgroups are 
 rational. A {\it rational polyhomomorphism} $R^\circ:G^\circ \rt H^\circ$ is a polyhomomorphism
 such that $\alpha(R^\circ)$, $\beta(R^\circ)$ are rational. 
Since indices in formulas (\ref{eq:product1})--(\ref{eq:product2}) are integer, 
 products of rational polyhomomorphisms are rational.
 
 \sm
 
 {\bf \punct  Polyendomorphisms.%
 \label{ss:polyend}}
 Now consider a unimodular locally compact group $G$ containing an open
 compact subgroup $K_0$. Normalize the Haar measure on $G$ assuming that the measure of $K_0$ is 1.
 For any pair $K_1\supset K_2$ of open compact subgroups 
 consider the index $[K_1:K_2]$. Consider the multiplicative semigroup $\Lambda=\Lambda(G)\subset \N$ consisting
 of all  products of such indices.
 Denote by $\Polh_\Lambda(G^\circ,G^\circ)$ the semigroup of all polyhomomorphisms
 $G^\circ\rt G^\circ$ consisting of 0 and
 all $R^\circ$ such that $\alpha(R^\circ)^{-1}$, $\beta(R^\circ)^{-1}\in \Lambda$.

 \begin{theorem}
 \label{th:rational} The set
  $\Polh_\Lambda(G^\circ,G^\circ)$ is a compact subsemigroup in $\Polh(G^\circ,G^\circ)$.
 \end{theorem}

The proof is contained in Subsect. \ref{ss:2rational}.

\sm

 {\bf \punct Example: the group of infinite matrices over a finite field.%
 \label{ss:example}}
 Let  $p$ be prime, $\F_p$ be the field with $p$ elements.
 Consider an infinite dimensional locally compact linear space over
  $\F_p$ satisfying the second county axiom.
  There are only 3 such spaces
  (this is semi-obvious, we present a formal proof below
  in Subsect.
 \ref{ss:remark}).
The first space 
 $\V_p^+$ is a direct sum of a finch number of copies of the field
   $\F_p$ equipped with a discrete topology.
   The second space  
 $\V_p^-$ is the direct product of a countable number of copies of the field
 $\F_p$, it is equipped with the Tikhonov topology. 
 These spaces are Pontryagin dual one to another
  (on the duality, see. e.g., \cite{HR}, Chapter 6).
  The third space  $\V_p:= \V_p^-\oplus \V_p^+$
  is the topic of our interest.

 Consider the linear space $\V_p$ over $\F_p$ consisting of two-side sequences
 \begin{equation}
 v=(\dots, v_{-2}, v_{-1}, v_0, v_1, v_2,\dots),\qquad v_k\in \F_p,
 \label{eq:v}
 \end{equation}
 such that $v_{j}=0$ for sufficiently large $j$. For each $m\in\Z$ consider the subspace
 $W^m\subset \V_p$ consisting of vectors $v$ such $v_l=0$ for all $l>- m$,
 $$
 \dots\supset W^{-1}\supset W^0\supset W^1\supset W^2\supset \dots.
 $$
  The topology
 in $\V_p$ is defined by the condition: the subgroups $W^m$ are open and form a base 
 of neighborhoods of 0. A sequence $v^{(l)}\in \V_p$
 converges to $v$ if it is contained in some subgroup $W^m$
 and converges to $v$ coordinate-wise.
 
  The subgroups $W^m$ are compact and are isomorphic
 to the countable direct product of cyclic groups $\Z_p$, quotients $\V_p/W^m$ are discrete
 and are isomorphic to the countable direct sum of of cyclic groups $\Z_p$.
 We normalize the Haar measure $\phi(v)$ on $\V_p$ assuming 
 that the measure of $W^0$ is 1. 
 
 \sm

 Consider groups 
$\GL(\V_p^+)$, $\GL(\V_p^-)$, $\GL(\V_p)$ of all continuous linear operators
in these spaces. The representation theory of $\GL(\V_p^+)$
is relatively simple (see T.~Tsankov \cite{Tsa}, it is not difficult
to reduce this classification to a result of G.~I.~Olshanski in \cite{Ols-semi}),  the group
$\GL(\V_p^-)$ is isomorphic to $\GL(\V_p^+)$. 
The group $\GL(\V_p)$ was introduced in \cite{Ner-finite}, the representation theory
of this group is non-trivial, see \cite{Ner-finite}, \cite{Ner-Wei}, \cite{Ner-GL},
it has many analogues with the representation theory of
infinite-dimensional real classical groups in the sense of
G.\,I.~Olshanski \cite{Ols-GB}.

 \sm

 Denote by   $\GL(\V_p)$ the group of all continuous linear operators in $\V_p$,
we can also say that it is  the group $\Aut(\V_p)$ of continuous automorphisms
 of the Abelian group $\V_p$.
 
 Denote by $J$ the operator of left shift  of sequences (\ref{eq:v}). Clearly,
 this transformation sends the Haar measure $\phi(v)$ to the measure $p\cdot \phi(v)$.
 Denote by $\GL^0(\V_p)$ the subgroup of $\GL(\V_p)$ consisting of transformations
 preserving the Haar measure on $\V_p$. Clearly, the group $\GL(\V_p)$
 is a semidirect product of the cyclic subgroup generated by $J$ and the normal subgroup
 $\GL^0(\V_p)$.

 \sm
 
We have a measure preserving action of $\GL^0(\V_p)$ by automorphisms of the locally compact group
$\V_p$, i.e., we are in the situation discussed in  Subsect. \ref{ss:polyend}.
The semigroup $\Lambda(\V_p)$ consists of powers $p^j$, where $j\ge 0$.
Closed subgroups in $\V_p\times\V_p$ are linear subspaces in $\V_p\oplus \V_p$.

\begin{theorem}
\label{th:GL}
 The closure of $\GL^0(\V_p)$ in $\Polh(\V_p,\V_p)$ coincides with
 the semigroup $\Polh_\Lambda(\V_p,\V_p)$.
\end{theorem}

The proof is contained in Subsect. \ref{ss:2GL}.

\begin{theorem}
	\label{th:last}
	Any unitary representation of the group $\GL^0(\V_p)$
	admits a continuous extension to a representation of the semigroup $\Polh_\Lambda(\V_p,\V_p)$
	compatible with the involution.
\end{theorem}	 
 
 The proof is contained in Subsect. \ref{ss:last}.

\sm 

{\bf \punct Problem of closure.}
Consider a unitary representation $\rho$ of a topological group $G$ in a Hilbert space $H$.
Consider  the set $\rho(G)$ of unitary operators and  close it in the space 
of all bounded operators with respect to the weak operator topology.
It can be readily checked that this closure $\ov{\rho(G)}$ 
is a compact semigroup. G.\,I.~Olshanski, see, e.g. \cite{Ols-semi}, showed that such semigroups
can be interesting algebraic objects and an effective tool for investigation
of unitary representations of infinite-dimensional groups $G$, see more in \cite{Ner-book}.

Now let a  group $G$ acts by  transformations of 
a measure space $X$. Then it acts in $L^2(X)$
and we have the same question about weak closure.
On the other hand such questions can be reformulated in the terms of
closures of groups in semigroups of polymorphisms,
apparently, the first problem of this type (closure of an infinite-dimensional orthogonal group
acting on a space of Gaussian measures) was solved by Nelson
\cite{Nel} (see, also \cite{Ner-Haar}, Sect. 12), for several actions of infinite-dimensional groups closures were described
in \cite{Ner-matching}, \cite{Ner-gauss}, \cite{Ner-Whi}. Theorem \ref{th:GL}
gives an additional example of this kind.

The problem of weak closure  is not interesting for semisimple real or $p$-adic groups
(usually, we get the one-point compactification, see  \cite{HM}).

On the other hand there are lot of interesting results about closures 
of ergodic measure preserving actions of Abelian groups as $\Z$, $\R$, $\Z^n$.
For a generic (in the sense of Baire category) transformations such closures
are huge and are related to the centralizer of a transformation in the semigroup
of Markov operators, see \cite{King}, \cite{JRR},  \cite{Sol}.
For non-mixing actions the problem of weak closure usually is a difficult
problems, see \cite{JPRV}, some relatively simple cases
 for spaces of infinite measures were examined in \cite{KR}, \cite{Ryzh}.

\sm 

{\sc Example.} 
a) Equip the countable space $\V_p^+$ (see the previous subsection) with the 
counting measure. It can be readily checked that the closure of $\GL(\V_p^+)$  
in $\Polh(\V_p^{+},\V_p^+)$ consists of partial linear bijections $\V_p^+\to \V_p^+$
equipped with counting measures. By \cite{Tsa}, the semigroup of partial linear 
bijections acts in all unitary representations of $\GL(\V_p^+)$.

\sm 

b) Equip the space $\V_p^-$ with the probabilistic Haar measure.
It is easy to show that the closure $\GL(\V_p^-)$ in $\Polh(\V_p^-,\V_p^-)$
consists of $R^\circ$ such that $\dom R=\V_p^-$, $\im R=\V_p^-$
and the Haar measure on $R$ is probabilistic.
\hfill $\boxtimes$

 \section{Proofs}
 
 \COUNTERS

 {\bf \punct Immediate corollories of the definition of polyhomomorphism. Proof of Proposition \ref{pr:add-def}.%
 \label{ss:prop1}} Here we prove statements a)-d) of Proposition  \ref{pr:add-def}.
 Recall that $R$ is a subgroup in $G\times H$ equipped with a left invariant Haar measure.
 
 \sm
 
 {\sc Statement a.} Let we have a locally compact group $N$,
 its closed normal subgroup $K$, and the quotient
 group $M$. Denote  left invariant Haar measures on these groups by
 $\nu(n)$, $\kappa(k)$, $\mu(m)$ respectively. 
 For $n\in N$ denote by $\dot n$ its image in $M$.
 According 
 \cite{Bour2}, Proposition VII.2.10, we have the following integration
 formula
 \begin{equation}
 \int_N f(n)\,d\nu(n)=\int_{M}\int_{K} f(\dot n k)\,d\kappa(k)\,d \mu(m).
 \label{eq:integration}
 \end{equation}
 Suppose that $K$ is not compact. Consider the image of the measure $\nu(n)$
 under the homomorphism $N\to M$. We wish to show that compact subsets $U$ in $M$ 
 with nonempty interiors have infinite measures. 
 Indeed, let $\wt U\subset N$ be the preimage of $U$. Applying the integration
 formula to the indicator function $I_{\wt U}$ we get $\infty$.
 
 We apply this remark to the group $N=R$, its subgroup $K=\indef R$, and the quotient
 $M=\dom R$ and observe that the projection of the measure
 $\rho(r)$ to $G$ can be dominated by the Haar measure  
 $\rho(g)$ only if $\indef R$ is compact. To verify the compactness of 
 $\ker R$, we take the group  $N=R$, the subgroup  $K=\ker R$ and the quotient group $M=\im R$.

\sm

{\sc Statement b. The image and the domain are open subgroups%
		\footnote{The author thanks the reviewer, who proposed a proof
			more natural than the proof in the first version
			of the paper. Initially, the statement was derived from the following 
			result of J.~Mackey \cite[Theorem 7.2]{Mack}: {\it  Let $G$ be a locally
				compact group, $L$ be a subgroup being a Borel subset
				in  $G$. Then the homogeneous space 
			$G/L$ has a standard Borel structure if and only if
			  $L$ is closed.} In our case the subgroup $\dom R$ has non zero measure,
		  therefore   $G/\dom R$ is at most countable. It easy to show that
		  the Borel structure on this set is standard.}.}
The group
 $R$ is a countable union of compact sets, 
 therefore its image 
$\dom R$ also is a countable union of compact sets, and therefore is a Borel set.

Let $A$ be a set of nonzero measure in a locally compact group.
Then by 
\cite{Wei}, \S11, the set  $AA^{-1}$ contains a neighborhood of the unit.
Hence a subgroup of nonzero measure contains a neighborhood of unit
and therefore is open.

\sm

{\sc Statement c. The group $R$ is unimodular.}  We apply formula (\ref{eq:integration})
to $N=R$, $K=\indef R$, $M=\dom R$.
 Since $\dom R$ is an open subgroup in a unimodular group $G$, it is unimodular. 
The subgroup $\indef R$ is compact, therefore it is unimodular,
moreover, the measure is preserved under all automorphisms of $\indef R$.
Formula (\ref{eq:integration}) shows that 
$\int f(hgh^{-1})=\int f(g)$ for all $h$, this implies that  $R$
is unimodular.

 Since
$R/\indef R\simeq \dom R$, we get that $R$ is unimodular.

\sm

{\sc Statement d.}
Let $g\in\dom R$. Then there is $h\in H$ such that $(g,h)\in R$.
The Haar measure in $R$ is invariant with respect to the left shift by $(g,h)$.
Therefore its projection $\nu$ to $G$ is invariant with respect to shift by $g$.
Therefore $\nu$ is a Haar measure on $\dom R$.

\sm

{\bf \punct Support of a product of polymorphisms.%
\label{ss:2support}} Here we derive Lemma \ref{l:proper}, which
is used below in discussion of  products of polyhomomorphisms.

Let $X$, $Y$ be compact separable metric spaces. We say that a relation
 $R:X\rra Y$ is {\it bi-proper}, if the projections 
 $R\to X$ are $R\to Y$ proper maps. 

An equivalent definition: a relation
 $R:X\rra Y$ is bi-proper if

\sm

1) $R\subset X\times Y$ is a closed subset;

\sm 

2) for any compact subset $A\subset X$ the set $RA\subset Y$  is compact;

\sm 

3) for any compact subset
 $B\subset Y$ the set $R^\square B\subset X$ is compact.

\begin{lemma}
	\label{l:proper-product}
	Let $R:X\rra Y$, $S:Y\rra Z$ be bi-proper relations. Then 
	$SR$ is bi-proper.
\end{lemma}

{\sc Proof.} It is sufficient to verify that  $SR\subset X\times Z$
is closed. Let a sequence
 $(x_j, z_j)\in X\times Z$ converges to  
 $(x^\circ,z^\circ)$. Then there is a sequence   $y_j\in Y$,
such that $(x_j,y_j)\in R$, $(y_j,z_j)\in S$. The set $\Xi\subset X$ consisting of
the sequence  $x_j$ and its limit  $x^\circ$ is compact.
Therefore the subset $R\Xi\subset Y$ is compact. 
This set contains the sequence
$(x_j,y_j)$, therefore we can choose a convergent subsequence    $(x_{j_k},y_{j_k})$.
Denote its limit by 
$(x^\circ,y^\circ)$. By the closeness of  $R$ this limit is contained  in $R$.  Clearly, $y_{j_k}$
converges to $y^\circ$, therefore $(y^\circ, z^\circ)\in S$. Hence 
$(x^\circ, z^\circ)\in SR$.
\hfill $\square$

\sm

Let 
 $X$, $Y$ be locally compact complete metric spaces
 equipped with measures    $\xi$, $\upsilon$ respectively. 
 We say that a polymorphism
$\mu:X\rt Y$ is  {\it bi-proper supported by}  
$R$, if

\sm

1. $R:X\rra Y$ is a bi-proper relation;  

\sm 

2. $\mu$ is supported by  $R$, i.e., $\mu\bigl((X\times Y)\setminus R\bigr)=0$.

\begin{lemma}
	\label{l:proper}
	Let $X$, $Y$, $Z$ be locally compact complete  metric spaces, 
	 $\xi$, $\upsilon$, $\zeta$ be measures on these subspaces.
	 Let a polymorphism
	$\mu\in\Pol(X,Y)$ be bi-proper supported by   $R$, 
	$\nu\in\Pol(Y,Z)$ be bi-proper supported by  $S$. 
	Then the product
	$\kappa:=\nu\mu$ is bi-proper supported by the product of relations $SR$. 
\end{lemma}

First, we prove the following lemma:

\begin{lemma}
	\label{l:support}
	Let $(X,\xi)$, $(Y,\upsilon)$ be locally compact complete metric  spaces with measures. 
	Let $\mu\in\Pol(X,Y)$ be a bi-proper polymorphism 
	supported by $R\subset X\times Y$. For a point
 $y_0\in Y$ consider the set  $R^\square y_0$, i.e. the set of all 
	$x\in X$ such that $(x,y_0)\in R$. Then for any neighborhood  $U$ 
	of the set
	$R^\square y_0$ there is a neighborhood  $V$ of the point $y_0$
	such that for any function 
	 $\phi$ supported by  $V$ the function $\Pi(\mu) \phi$
	has a support in $U$.
\end{lemma}	

{\sc Proof of Lemma \ref{l:support}.}
Let find
$\Pi(\mu)\phi$ applying Remark from Subsect.  \ref{ss:polymorphisms}.B.
For this purpose we must project the measure
 $\psi(y) \mu(x,y)$ to the space
 $X$. Clearly, for a function  $\phi$ whose support is contained in a small
 neighborhood of $y_0$, the support of the measure $\psi(y) \mu(x,y)$
 is contained in a small neighborhood of the set
 $(X\times y_0)\cap R$. The support of the projection of the measure $\psi(y) \mu(x,y)$
 is contained in a small neighborhood of the set
 $R^\square y$.
\hfill $\square$

\sm 

{\sc Proof of Lemma \ref{l:proper}.} 
Let $(x_0,z_0)\notin SR$. We must show that for real functions 
$\phi\in C_c(X)$ supported by a  sufficiently small neighborhood 
$A$ of the point
$x_0$ and $\theta\in C_c(Z)$ supported by a  sufficiently small neighborhood 
  $B$ of the point 
$z$, we have $\int_{X\times Z} \phi(x)\,\theta(z)\,d\kappa(x,z)=0$.
Evaluating this expression, we get
\begin{multline}
\int_{X\times Z} \phi(x)\theta(z)\,d\kappa(x,z)=
\bigl\la \phi,\Pi(\kappa)\theta\bigr\ra_{L^2(X,\xi)}
=
\bigl\la \phi,\Pi(\mu)\Pi(\nu)\theta\bigr\ra_{L^2(X,\xi)}
=\\=
\bigl\la \Pi(\mu^\square)\phi,\Pi(\nu)\theta\bigr\ra_{L^2(Y,\upsilon)}.
\label{eq:phitheta}
\end{multline}
By Lemma
 \ref{l:support} a support of the function  $\Pi(\nu)\theta$
 is contained in a small neighborhood of the set
 $S^\square z_0$,
 a support of the function
 $\Pi(\mu^\square)\phi$ is contained in a small neighborhood
 of the set  $R x_0$. The condition $(x_0,z_0)\notin SR$
 is equivalent to
 $S^\square z_0\cap R x_0=\varnothing$.
 Therefore \eqref{eq:phitheta} is zero for functions $\phi$, $\theta$
 with sufficiently small supports.

Other statements follow from Lemma
 \ref{l:proper-product}.
\hfill $\square$

\sm

{\bf \punct The product of polyhomomorphisms.%
\label{ss:2products-poly}}
Thus, we have two polyhomomorphisms $R^\circ:G^\circ\rt H^\circ$, $T:H^\circ\to K^\circ$.
Notice that the indices
$$[\ker T:(\ker T\cap \im R)], \qquad \indef R:(\indef R\cap \dom T)$$
in formulas (\ref{eq:product1})--(\ref{eq:product2}) is finite. Indeed, the subgroup $\im R$
is open (and therefore closed). Hence $\im R\cap \indef T$
is open and closed in the compact group $\indef T$. Therefore the 
quotient $(\indef T)/(\im R\cap \indef T)$ is finite.

\sm

Next, {\it we intend to  evaluate the  product of $R^\circ$ and $T^\circ$ as a product of polymorphisms
and to verify that it coincides with the product $T^\circ R^\circ$ in the sense of polyhomomorphisms.}
By Lemma \ref{l:proper} the product is supported by
the closed subgroup $RT$. We must show that the measure on $RT$ is a Haar measure.

Denote by $L_G(u)$ the transformation $v\mapsto uv$ on a group $G$.
For the polyhomomorphism $R^\circ$ we have the following identity for polymorphisms:
$$L_H(h) R^\circ= R^\circ L_G(g), \qquad \text{for $(g,h)\in R$}.
$$
If also $(h,k)\in T$, then we have
$
T^\circ L_H(h) = L_K(k) T^\circ
$
and therefore 
$$ R^\circ T^\circ L_G(g)=L_K(k) R^\circ T^\circ .$$
Therefore the polymorphism $R^\circ T^\circ$
is determined by a Haar measure
on $RT$. It remains to find normalization constants
$\alpha(T^\circ R^\circ )$, $\beta(T^\circ R^\circ )$.

\begin{lemma}
	\label{l:I}
Let $R^\circ\in\Polh(G^\circ, H^\circ)$, let $\Pi(R^\circ):L^2(H,\eta)\to L^2(G,\gamma)$ be the
corresponding operator.

\sm

{\rm a)} Let $B\subset \im R$ be a compact subset of nonzero measure invariant with respect
to $\indef R$. Then
$$
\Pi(R^\circ) I_B=\alpha(R^\circ) I_{R^\square B}.
$$


{\rm b)} Let $Z\subset \indef R$ be a subgroup of finite
index  $N$.
Let $h_1$, \dots, $h_N$ be representatives of double cosets --  $\indef R/Z$.
Let $D\subset \im R$ be a compact set of non-zero measure 
invariant with respect to  $Z$. Let the sets $h_j C$  are mutually disjoint. Then 
$$
\Pi(R^\circ) I_D=
\frac{\alpha(R^\circ)}N I_{R^\square D}.
$$

\end{lemma}

{\sc Remarks.} a) Since the subgroup $\indef R$ is normal in  $\im R$, the left $(\indef R)$-invariance
of   $B$ is equivalent to 
right   $(\indef R)$-invariance.

\sm 

b) For any compact $D$ we have
\begin{equation}
\qquad\qquad\qquad\Pi(R^\circ) I_D= \Pi(R^\circ) I_{D\cap \im R}.\qquad\qquad\qquad
\boxtimes
\label{eq:ID}
\end{equation}

{\sc Proof.} a) We evaluate $\Pi(R^\circ) I_B$ using Remark in Subsect. \ref{ss:polymorphisms}.b.
The measure (\ref{eq:psimu}) is supported by the set $M$ of all $(a,b)\in G\times H$ such that 
$(a,b)\in R$, $b\in B$ and coincides with the Haar measure $\rho(r)$ on this set.
If $(a,b)\in M$ and $q\in \indef R^\circ$, then $(a,bq)\in M$.
Therefore $M$  coincides with the preimage of $R^\square B$ under the 
projection $R\to G$. Projecting the Haar measure on $M$ to $G$
we get the measure $\alpha(R^\circ)\,\gamma$ restricted to
$R^\square B$. 

\sm

b) Obviously, for any set  $D\subset H$ and $h\in \indef R$ 
the followin equality holds
$$\Pi(R^\circ)I_{hD}=\Pi(R^\circ) I_D.$$
Therefore, in our case all functions
 $ \Pi(R^\circ)I_{h_jC}$ are equal. 
The sets  $R^\square (h_j C)$ also coincide. Therefore
$$
\Pi(R^\circ)I_{C}=\frac 1N\sum_j \Pi(R^\circ)I_{h_jC}=\frac 1N\sum_j \Pi(R^\circ)I_{\cup h_jC}.
$$
Applying the statement a) to the last set, we come to
$$
\qquad\qquad\qquad
\frac {\alpha(R^\circ)}N I_{R^\square (\cup h_jC)}=
\frac {\alpha(R^\circ)}N I_{\cup R^\square  h_jC}=\frac {\alpha(R^\circ)}N I_{R^\square C}.\qquad\qquad\quad
\square
$$


\sm

{\sc End of proof of Theorem \ref{th:product}.}
Consider a compact subset 
$C\subset  \im (TR)$ containing a neighbothood
of the unit and invariant with respect to the subgroup
$$\indef (TR)= T \indef R.$$
By virtue of Lemma
\ref{l:I}.a,
$$\Pi((TR)^\circ) I_C
=\alpha((TR)^\circ)I_{R^\square T^\square C}.
$$
By the same lemma and
 \eqref{eq:ID}, we have
$$
\Pi(R^\circ)\Pi(T^\circ) I_C=\alpha(T^\circ)\Pi(R^\circ) I_{T^\square C}=\alpha(T^\circ)\Pi(R^\circ) I_{\im R\cap T^\square C}.
$$
Generally, the set
 $\im R\cap T^\square C$ is not
$(\indef R)$-invariant. However, by our choice of
 $C$, it is $(\indef R\cap \dom T)$- invariant.
 Applying Lemma
 \ref{l:I}.b to the last expression, we come to
$$
\frac{\alpha(T^\circ)\alpha(R^\circ)}{[\indef R:(\indef R\cap \dom T)} I_{R^\square T^\square C}.
$$

The same considerations give us the constant $\beta(R^\circ T^\circ)$.
\hfill $\square$

\sm

{\bf\punct The description of operators $\Pi(R^\circ)$.%
\label{ss:2operators}}
Let $G$ be a locally compact group, $\Phi$ an open subgroup, and $\Delta$ a  compact normal
subgroup in
$\Phi$. Normalize a Haar measure on the group $\Phi/\Delta$ as the image of $\gamma(g)\bigr|_{\Phi}$
under the map $\Phi\to\Phi/\Delta$. Consider the 'diagonal' map $\Phi\to \Phi\times (\Phi/\Delta)$
sending $g\in\Phi$ to $(g,g\Delta)$. We define the polyhomomorphism
$$\mu^\circ_G[\Phi|\Delta]\in\Polh\bigl(G^\circ,(\Phi/\Delta)^\circ\bigr)$$
as the image of $\gamma(g)\bigr|_{\Phi}$ under the 'diagonal' map
(in particular $\alpha=\beta=1$).

On the other hand the map $\Phi\to \Phi/\Delta$
induces the operator
$$\Pi(\mu^\circ_G[\Phi|\Delta]): L^2(\Phi/\Delta)\to L^2(\Phi)\subset L^2(G).$$
It is an isometric embedding $L^2(\Phi/\Delta)\to L^2(G)$
whose image is $L^2(\Phi)^\Delta$. The adjoint operator
$$
\Pi(\mu^\circ_G[\Phi|\Delta]^\square):L^2(G)\to L^2(\Phi/\Delta)
$$
can be described in the following way: we restrict a function
$f\in L^2(G)$ to the open subgroup $\Phi$, take its average over 
the action of the compact group $\Delta$, and consider this average as a function
on $\Phi/\Delta$. 

\sm 

Let $R^\circ\in \Polh[G^\circ,H^\circ]$. We decompose it as a product
$R^\circ=T^\circ S^\circ Q^\circ$, 
$$
G^\circ\, \overset{Q^\circ}\rt \,(\dom R/\ker R)^\circ\, \overset{S^\circ}\rt\,
(\im R/\indef R)^\circ\, \overset{T^\circ}\rt\, H^\circ
, 
$$
where
 $$
 Q^\circ=\mu^\circ_G[\dom R|\ker R],\qquad
T^\circ=\mu^\circ_H[\im R|\indef R]^\square.
$$
To define $S^\circ$ we consider the  canonical  map 
$$R\to (\dom R/\ker R)\times (\im R/\indef R).$$
Its image is a graph of an isomorphism 
$$\Sigma:(\dom R/\ker R)\to (\im R/\indef R).$$
The Haar measure on this graph is the image of the measure $\rho(r)$.

The operators $\Pi(Q^\circ)$,  $\Pi(T^\circ)$ were described above in this subsection.
$$
\Pi(S^\circ) f(q)=\beta(R^\circ)\cdot f(\Sigma(q)), \qquad \text{where $q\in \im R/\indef R.$}
$$






{\sc Proof of Proposition \ref{pr:homothety}.} Thus, 
the operator   $\Pi(R^\circ)$ decomposes
into the product of three operators. The operators $\Pi(T^\circ)$
is the operator of projection 
\begin{multline*}
L^2(H)\simeq L^2(\im R)^{\indef R}\oplus
\bigl(L^2(\im R)^{\indef R} \bigr)^\bot
\to \\\to
L^2(\im R)^{\indef R}\simeq L^2(\im R/\indef R)
\end{multline*}
Next, we apply the operator
$$\Pi(S^\circ): L^2(\im R/\indef R)\to L^2(\dom R/\ker R),$$
which is unitary up to a scalar factor.
The last operator is an isometric embedding
\begin{multline*}
L^2(\dom R/\ker R)\simeq L^2(\dom R)^{\ker R}
\to\\\to L^2(G)\simeq
L^2(\dom R)^{\ker R}\oplus \bigr(L^2(\dom R)^{\ker R} \bigl)^\bot.
\end{multline*}
Clearly, the product is a partial isometry up to a scalar factor,
the initial subspace is
$L^2(\im R)^{\indef R}$, the final subspace is $L^2(\dom R)^{\ker R}$.

\sm

{\bf\punct The proof of Proposition  \ref{pr:angles}.
\label{ss:angles-proof}}
 Recall that  $\Phi$, $\Psi$ are open subgroups in  $G$,
$\Delta$ is a normal subgroup in $\Phi$, and $\Gamma$ is a normal subgroup
in $\Psi$.
Set
$$V:= L^2(\Phi)^\Delta \quad \text{and}\quad W:=L^2(\Psi)^\Gamma\subset L^2(G).$$
Denote by
$$P:=P_{\Phi|\Delta}^G,\qquad Q:=P_{\Psi|\Gamma}^G$$
the projection
operators to these subspaces.
{\it We wish to show that the self-adjoint operator
	$$PQ\Bigr|_V=PQP\Bigr|_V:V\to V$$
	splits as a direct sum of a zero operator and a scalar operator.}
We can pass to the subspaces
$V\ominus (V\cap W^\bot)$, $W\ominus (W\cap V^\bot)$.
Indeed, $Q$ is zero on $V\cap W^\bot$, and  $W\cap V^\bot$
does not contained in the image of
 $Q$. Therefore, it is sufficient to show that
$
PQ\Bigr|_{V\ominus (V\cap W^\bot)}
$
is a scalar operator.

Next, let a function
$f\in V$ has a support in  $\Phi\setminus \Psi$. Obviously, 
$Qf=0$. But  $f\in V$ is  $\Delta$-invariant, therefore 
a support of $f$ is actually contained in 
$$\Phi\setminus \Delta(\Phi\cap\Psi)=\Phi\setminus (\Phi\cap\Psi)\Delta.$$
{\it Therefore, without loss of a generality, we can assume that}
$$
\Phi=\Delta\cdot (\Phi\cap \Psi),\qquad\Psi=\Gamma\cdot (\Phi\cap \Psi).
$$
Notice that (under this condition)
$f\in V$ is determined by its value on 
$\Phi\cap \Psi$, and this restriction is $\Delta\cap \Psi$-invariant.
Conversely, any $\Delta\cap \Psi$-invariant function from  $L^2(\Phi\cap\Psi)$
admits an extension to the whole subgroup
 $\Phi$ by $\Delta$-invariance
(and by 0 on $G\setminus \Phi$).
By the way, for two functions  $f_1$, $f_2\in V$ the following identity holds
$$
\la f_1,f_2\ra_{L^2(G)}=\Bigl\la f_1\Bigr|_\Phi,f_2\Bigr|_\Phi\Bigr\ra_{L^2(\Phi)}=[\Delta:(\Delta\cap \Psi)] \cdot \Bigl\la f_1\Bigr|_{\Phi\cap \Psi}, f_2\Bigr|_{\Phi\cap \Psi}\Bigr\ra_{L^2(\Phi\cap \Psi)}.
$$
Similar statements take place for elements of the subspace
 $W$. 

If  $f\in V$, $h\in W$, then
$$
\la f,h\ra_{L^2(G)}= \Bigl\la f\Bigr|_{\Phi\cap \Psi}, h\Bigr|_{\Phi\cap \Psi}\Bigr\ra_{L^2(\Phi\cap \Psi)}
$$
The operator
 $Q\Bigr|_V$ is described in the following way: we take a function
 $f\in V$ and average with respect to  $\Gamma$ the function 
$f\Bigr|_{\Phi\cap \Psi}$.
Point out that
\begin{equation}
Qf\Bigr|_{\Phi\cap \Psi}=\frac 1{[\Gamma:(\Gamma\cap \Phi)]}
\Bigl\{\text{average of $f\Bigr|_{\Phi\cap \Psi}$ with respect to $\Gamma\cap \Phi$} \Bigr\}.
\label{eq:average}
\end{equation}
Notice, that the subgroups
 $\Delta\cap \Psi$ and $\Gamma\cap \Phi$
are normal in $\Phi\cap \Psi$, therefore the averaging with respect to  $\Gamma\cap \Phi$
sends $\Delta\cap \Psi$-invariant functions to $\Delta\cap \Psi$-invariant functions.

\begin{lemma}
	Let $f\in V\ominus (V\cap W^\bot)$. Then the function  $f\bigr|_{\Phi\cap\Psi}$
is invariant with respect to the subgroup	
	$(\Delta\cap \Psi)\cdot (\Gamma\cap \Phi)$.  
\end{lemma}	

{\sc Proof.} 
The space $V\cap W^\bot$ consists of functions $f$ such that
the expression~\eqref{eq:average} is 0. Therefore $V\ominus(V\cap W^\bot)$
consists of $\Gamma\cap \Phi$-invariant elements of  $V$. \hfill $\square$

\sm

\sm

{\sc Proof of Proposition  \ref{pr:angles}.}
For a function $f\in V\ominus (V\cap W^\bot)$ denote by  $\wt f$ its restriction to  $\Phi\cap\Psi$.
Formula~\eqref{eq:average} gives
$$Qf\Bigr|_{\Phi\cap \Psi}=[\Gamma:(\Gamma\cap \Phi)]^{-1} \wt f.$$
The similar formula for the operator
 $P$ leads to 
$$
Pf\Bigr|_{\Phi\cap \Psi}=[\Delta:(\Delta\cap \Psi)]^{-1}\, [\Gamma:(\Gamma\cap \Phi)]^{-1} \wt f,
$$
and we get the desired statement.
\hfill $\square$

\sm


{\bf\punct Proof of Theorem \ref{th:rational} (semigroups $\Polh_\Lambda(G^\circ,G^\circ)$ are closed).%
	\label{ss:2rational}}
	Let a sequence of polyhomomorphisms
	 $R_j^\circ=(R_j,\rho_j)$ containing in the semigroup  $\Polh_\Lambda(G^\circ,G^\circ)$
	converges to  $R^\circ=(R,\rho) \in \Polh(G^\circ,G^\circ)$.
	We must show that $R^\circ \in \Polh_\Lambda(G^\circ,G^\circ)$.

	By symbols
	$\rho_j$, $\rho$ we also denote the corresponding measures
	on  $G\times G$. 
Let $R_j^\circ\in \Polh_\Lambda(G^\circ,G^\circ)$ converge to $R^\circ \in \Polh(G^\circ,G^\circ)$.
Without loss of generality we can assume that the sequence $\alpha(R^\circ_j)$
converges (otherwise we pass to a subsequence). If it converges to 0, then $R^\circ_j$ converges to the zero polyhomomorphism.
Otherwise, $\alpha(R^\circ_j)$ is eventually constant, without loss of generality we can assume that 
it is constant. Also, we can assume that a sequence $\beta(R^\circ_j)$ is constant.

Consider a compact open subgroup  $L\subset \dom R^\circ$ containing $\ker R$.
Denote $M:=RL$. Then $L\times M$ is a compact open subgroup in $G\times H$.
We have
$$\rho\Bigr|_{L\times M}=\lim_{j\to\infty} \rho_j\Bigr|_{L\times M}$$
Next $\pil (R_j\cap (L\times M))$ is an open subgroup
in $L$ of a certain index $p_j$. Without loss of generality we can assume that
$p_j$ does not depend on $j$ or $p_k\to \infty$. 
In the second case we have 
$$
\gamma\Bigl(\pil (\rho_j\bigr|_{L\times M}\Bigr)\le \gamma\Bigl(\pil (\rho_j\bigr|_{L\times H}\Bigr) =
\frac \alpha{p_j}\gamma(L)
$$
and $R^\circ_j$ converges to zero. So we consider the first case.

Groups $R_j\cap (L\times M)$ are open subgroups in $R_j\cap (L\times H)$
of indices $q_j$. Again, without loss of generality, we can assume that
the sequence 
$q_j$ is constant. Now we have
$$
\gamma\Bigl(\pil \bigl(\rho_j\bigr|_{L\times M}\bigr) \Bigr)=\frac{\alpha}{pq}\, \gamma(L).
$$
Therefore, $\rho_j(L\times M)=\frac{\alpha}{pq} \gamma(L)$.
Passing to the limit, we get $\rho(L\times M)=\frac{\alpha}{pq} \gamma(L)$
and $\alpha(R^\circ)=\frac 1{pq} \alpha$.
\hfill $\square$

\sm

{\bf\punct Locally compact linear spaces over  finite fields.\label{ss:remark}}
Here we prove the statement about classification of locally compact linear spaces
over the field $\F_p$, which was formulated at the beginning 
of Subsect.~\ref{ss:example}.

Thus, let  $V$ be an infinite locally compact linear space over  $\F_p$,
equivalently, 
$V$ is a locally compact group such that
$p\cdot v=0$ for all $v\in V$.
Three cases are possible. 

\sm

{\it	The first case.} Let  $V$ be discrete   and therefore  countable.
Any countable linear space over
$\F_p$ is isomorphic to  $\V_p^+$.

\sm

{\it The second case.} 
Let $V$ be compact and infinite. Since $V$ is compact, the Pontryagin dual
group  $V^\circ$ is discrete. Therefore 
$V$ is dual to $\V_p^+$, i.e., $V\simeq \V_p^-$.

\sm 

{\it The third case.}
Let $V$ by non-compact and non-discrete. 
For any character $\chi$ from $V$
to the multiplicative groups of complex numbers the identity $\chi(v)^p=1$
holds, i.e., values of    $\chi$
have the form $e^{2\pi i/p}$. According the Pontryagin duality,
characters separate points of  $V$. Hence   $V$
is totally disconnected and therefore 
 $V$ contains an open compact subgroup 
$W$ (see \cite[Theorem~7.5]{HR}). If  $W$ is finite, then  $V$ is countable.
If   $V/W$ is finite, then  $V$ is compact.
So we can omit these cases. Thus, $W\simeq \V_p^-$, $V/W\simeq \V_p^+$.
Next, we choose a basis
$e_j$ in $V/W$ and choose representatives
$\wt e_j\in V$. Then the linear span
of vectors
 $\wt e_j$ is a discrete linear space 
 in  
 $V$ complementary to $W$.

\sm

{\bf \punct The closure of $\GL^0(\V_p)$ in the semigroup of polyhomomorphisms.%
\label{ss:2GL}} Here we prove Theorem \ref{th:GL} about the closure 
of the group $\GL^0(\V_p)$ in the semigroup of polyhomomorphisms of the space $\V_p$.

We must   show that each element of $\Polh_\Lambda(\V_p,\V_p)$
is contained in the closure $\ov{\GL}$ of the group $\GL^0(\V_p)$.
Point out that proofs in the present and  next subsections
are essentially based on results of
 \cite{Ner-GL}.

For $m>0$
denote by $\theta_m:\V_p\rt \V_p$ the  linear relation
consisting of $(v,v')$ such that
$v_j=v_j'=0$ for $j\ge m$, 
$v_i=v_i'$ for $-m<i<m$, and
$v_j$, $v_j'$ are arbitrary if $v_j\le -m$.
Thus,
$$
\ker \theta_m=W^m,\quad \indef \theta_m=W^m,\quad
\dom \theta_m=W^{-m},\quad \im \theta_m=W^{-m},
$$
and 
$$
\dom \theta_m/ \ker \theta_m\simeq \F_p^{2m-1}\simeq \im \theta_m/\indef \theta_m,
$$
where $\F_p^{2m-1}$ consists of vectors $(v_{-m+1}, v_{-m+2},\dots, v_{m-1})$.
The isomorphism
$\dom \theta_m/ \ker \theta_m\to  \im \theta_m/\indef \theta_m$
is the identity map $\F_p^{2m-1}\to \F_p^{2m-1}$.

We define the  polyhomomorphisms $\theta^\circ_m\in \Polh(\V_p,\V_p)$
assuming that $\alpha(\theta^\circ_m)=\beta(\theta^\circ_m)=1$.

\begin{lemma}
	{\rm a)} $\theta_m^\circ \in \ov\GL$.
	
	\sm 
	
	{\rm b)} The sequence $\theta_m^\circ$ converges to the identical
	polyhomomorphism as $m$ tends to $\infty$.
\end{lemma}

{\sc Proof.} Decompose the space $\V_p$ as a product of 3 measure spaces
$$\V_p=V^-\times \F^{2m-1} \times V^+,$$
where $V^-$ consists of sequences $(\dots, v_{m-1},v_m)$,
the $\F_p^{2m-1}$ of vectors 
$$(v_{-m+1}, v_{-m+2},\dots, v_{m-1}),$$
and $V^+$ of vectors $(v_m, v_{m+1}, \dots)$.
The space $V^+$ is countable and measures of all points
are $1$. The space $\F_p^{2m-1}$ is finite and  measures
of all points are $p^{-m}$. The space $V^-$ is equipped by
product of uniform probabilistic measures on $\F_p$. The measure of $W^0$ is 1.

Consider the sequence
$$S_j^+=
\begin{pmatrix}
0&1_j&0\\
1_j&0&0\\
0&0&1_\infty
\end{pmatrix}
$$
of linear transformations in $V^+$.
Clearly, it converges in $\Polh(V^+,V^+)$
to the delta measure supported by $0$.

Next, consider the sequence 
$$S_j^-:=\begin{pmatrix}
1_\infty&0&0\\
0&0&1_j\\
0&1_j&0
\end{pmatrix} 
$$
of linear transformations of the space
$V^-$. Clearly, in $\Polh(V^-,V^-)$ it converges to the
product measure on $V^-\times V^-$.

Let us regard $S_j^+$ (resp. $S_j^-$) as polymorphisms of the whole $\V_p$.
Then
$$
\lim_{j\to\infty}\lim_{i\to\infty} S^+_i S^-_j =\theta_m^\circ.
$$

b) Consider the compact subgroup $W^l$ with $l>0$,
a vector $v\in W^{-k}$, and the indicator function
$I_{v+W^l}$. Then for $m>\max(l,k)$ we have $\Pi(\theta^\circ_m)I_{v+W^l}=I_{v+W^l}$.
\hfill$\square$

\sm 

By the separate continuity of the product the statement a) implies the following corollary.

\begin{lemma}
	\label{l:1}
For any $g\in \GL^0$ we have $\theta^\circ_m g \theta^\circ_m\in \ov{\GL}$.
	\end{lemma}

\begin{lemma}
	\label{l:2}
Let $R^\circ\in \Polh_\Lambda(\V_p,\V_p)$.
Then there exists $g\in \GL^0$ such that  
$$\theta^\circ_m R^\circ \theta^\circ_m=\theta^\circ_m g \theta^\circ_m.$$	
Moreover, we can choose a finitary $g$, i.e., $g$ such
that $g-1$ has only finite number of nonzero matrix elements.  	
\end{lemma}

{\sc Proof.} Notice that for $Q^\circ_m:=\theta^\circ_m R^\circ \theta^\circ_m$,
$$
\dom Q_m\subset W^{-m},\quad \im Q_m\subset W^{-m}, \quad
\ker Q_m \supset W^m, \quad \indef Q_m\supset W^m. 
$$
Therefore $Q^\circ_m$ determines a polyhomomorphism
$$
W^{-m}/W^m\rt W^{-m}/W^m,
$$
i.e, $\F^{2m-1}\rt \F^{2m-1}$, measures on both copies of $\F^{2m-1}$
are uniform, a measure of a point is $p^{-m}$.

In particular, we can apply this reasoning to $\theta^\circ_m g \theta^\circ_m$,
where $g$ is finitary matrix. This polyhomomorphism determines
a polyhomomorphism $\chi^\circ(g):\F^{2m-1}\rt \F^{2m-1}$.
The  underlying linear relation $\chi(g):\F^{2m-1}\rra \F^{2m-1}$
 consists of $(u,v)$ 
such that there exist $x\in W^m$, $y\in W^m$ satisfying
$$
\begin{pmatrix}
x\\u\\0
\end{pmatrix}
\begin{pmatrix}
g_{11}&g_{12}&g_{13}\\
g_{21}&g_{22}&g_{23}\\
g_{31}&g_{32}&g_{33}
\end{pmatrix}
\begin{pmatrix}
y\\v\\0
\end{pmatrix}.
$$
This means that  $\chi(g)$ is the {\it characteristic linear relation}
of $g$ in the sense \cite{Ner-GL}, Subsect. 1.5.
Next, we must find  the normalization of the Haar measure $\theta^\circ_m g \theta^\circ_m$.
Evaluating $\alpha(\theta^\circ_m \cdot g \theta^\circ_m)$ by formula (\ref{eq:product1})
we get 
$$\alpha(\chi^\circ(g))=p^{-\rk g_{13}}.$$
In notation of \cite{Ner-GL}, Subsect. 1.5, $\rk g_{13}$
is the invariant $\eta(g)$.
So we get a polyhomomorphism $\chi^\circ(g):\F^{2p-1}\rt \F^{2p-1}$
such that measure of each point of $\chi(g)$ is $p^{-m-\rk g_{13}-\dim\indef \chi(g)}$
and 
$$\beta(\chi^\circ(g))=p^{-\rk g_{13}-\dim\indef \chi(g)+\dim \ker(g)}.$$

Now consider an arbitrary linear relation $Q:\F^{2m-1}\rra \F^{2m-1}$
and a polyhomomorphism $Q^\circ$ with $\alpha(Q^\circ)=p^{-\mu}$.
Then 
$$\beta(Q^\circ)=\alpha(Q^\circ)\,p^{\dim\ker Q-\dim\indef P}$$
We have $\beta\le 1$. By \cite{Ner-GL}, Proposition 1.8, any such polyhomomorphism 
can arise as $\chi^\circ(g)$ for a finitary $g$.
\hfill $\square$

\begin{lemma}
	\label{l:3}
For any $R^\circ \in \Polh(\V_p,\V_p)$ the sequence
$Q_m^\circ:=\theta_m^\circ R^\circ	\theta_m^\circ$ converges to $R^\circ$.
Also
$$
\lim_{n\to\infty, \,m\to\infty} \theta_m^\circ R^\circ	\theta_n^\circ=R^\circ.
$$
\end{lemma}	

{\sc Proof.} Fix $W_k$ and two vectors $v$, $w\in W_{-l}$.
Clearly the sequence
$$Q_m^\circ\bigl((v+W_k)\times (w+W_k)\bigr)$$
became constant after $m=\max(k,l)$.
\hfill $\square$ 

\sm 

Theorem \ref{th:GL} follows from Lemmas
\ref{l:1}--\ref{l:3}.

\sm 

{\bf \punct Semigroup extensions of unitary representations of $\GL^0(\V_p)$.%
\label{ss:last}}
Below we give the proof of Theorem \ref{th:last}, it is based
on \cite{Ner-book}, Theorem VIII.1.10.

We define the category $\ov \cK$, whose objects are spaces
$\F^{2m-1}_p$ with Haar measure normalized as above and $\V_p$. Morphisms are polyhomomorphisms.
We define the subcategory $\cK$, whose objects are the spaces $\F^{2m-1}_p$
with the same morphisms.

For any $m<n<\infty$ we define the linear relation $\lambda_{mn}:\F^{2m-1}\rra \F^{2n-1}$
as the subspace consisting of vectors
$$
(v_{-m+1}, \dots, v_{m-1})\oplus (v_{-n+1},\dots, v_{-m}, v_{-m+1}, \dots, v_{m-1},0,\dots,0).
$$
We also define linear relations $\lambda_{m\infty}:\F^{2m-1}\rra \V_p$
consisting of vectors 
$$
(v_{-m+1}, \dots, v_{m-1})\oplus (\dots, v_{-m-1}, v_{-m}, v_{-m+1}, \dots, v_{m-1},0, 0\dots).
$$
Define corresponding polyhomomorphisms $\lambda_{mn}^\circ$, $\lambda_{m\infty}^\circ$
assuming that all $\alpha(\cdot)$, $\beta(\cdot)$ are 1.
Define adjoint polyhomomorphisms $\mu_{mn}^\circ:=(\lambda_{mn}^\circ)^\square$, 
$\mu_{m\infty}^\circ=(\lambda_{m\infty}^\circ)^\square$.
It can be readily checked that we get a structure of an
 ordered category in the sense of \cite{Ner-book}, Sect. III.4.
 
 By \cite{Ner-GL} any unitary representation of $\GL^0(\V_p)$
 generates a representation of the category $\cK$.
 Our
Lemma \ref{l:3} allows to apply Approximation theorem
VIII.1.10 from \cite{Ner-book}, this implies that any $*$-representation
of the category $\cK$ extends to a representation of $\ov \cK$.
In particular the representation of the group $\GL^0(\V_p)$ extends
to a representation of $\Polh_\Lambda(\V_p,\V_p)$.

 \tt
\noindent
Yury Neretin\\
Wolfgang  Pauli Institute/c.o. Math. Dept., University of Vienna \\
\&Institute for Theoretical and Experimental Physics (Moscow); \\
\&MechMath Dept., Moscow State University;\\
\&Institute for Information Transmission Problems;\\
yurii.neretin@math.univie.ac.at
\\
URL: http://mat.univie.ac.at/$\sim$neretin/


\begin{thebibliography}{cc}

\bibitem{Bir}
Biringer, I.
{\it Metrizing the Chabauty topology.}
Geom. Dedicata 195 (2018), 19-22.

\bibitem{Bog}
Bogachev V. I. {\it Measure theory.} Vol. II. Springer-Verlag, Berlin, 2007.


\bibitem{Bour1}
Bourbaki N. {\it General topology. Chapters 1-4.} 
Reprint of the 1989 English translation. Elements of Mathematics (Berlin). Springer-Verlag, Berlin, 1998

\bibitem{Bour2}
Bourbaki N. {\it Integration. II. Chapters 7-9.} Translated from the 1963 and 1969.
Elements of Mathematics (Berlin). Springer-Verlag, Berlin, 2004. 

\bibitem{JPRV}
Janvresse, \'E.; Prikhod'ko, A. A.; de la Rue, T.; Ryzhikov, V. V. {\it Weak limits of powers of Chacon's automorphism.} Ergodic Theory Dynam. Systems 35 (2015), no. 1, 128-141.

\bibitem{JRR}
Janvresse \'E, de la Rue Th., Ryzhikov V.
{\it Around King's rank-one theorems: flows and $\Z_n$-actions.} in {\it Dynamical systems and group actions}
(eds.  Bowen L.,  Grigorchuk R.,  Vorobets Y.), 143-161,
Contemp. Math., 567, Amer. Math. Soc., Providence, RI, 2012.

\bibitem{HR}
Hewitt E., Ross K. A. {\it Abstract harmonic analysis. Vol. I: Structure of topological groups. Integration theory, group representations.}
Academic Press, New York, 1963 

\bibitem{Hopf}
 Hopf, E.
{\it The general temporally discrete Markoff process.} J. Ration. Mech. Anal.,
 3 (1954)
 13-45.

\bibitem{HM}
Howe R. E., Moore C. C.
{\it Asymptotic properties of unitary representations.} 
J. Funct. Anal. 32, 72-96 (1979).

\bibitem{Kech}
Kechris A. S. {\it Classical descriptive set theory.}
 Springer-Verlag, New York, 1995. 
 
 \bibitem{King}
 King J.
{\it The commutant is the weak closure of the powers, for rank-1 transformations}.
Ergodic Theory Dynam. Systems 6 (1986), no. 3, 363-384.

 \bibitem{Kre}
 Krengel U. {\it Ergodic theorems.} Walter de Gruyter \& Co., Berlin, 1985.
 
 \bibitem{KR}
 Kushnir A. Yu.,  Ryzhikov V. V., {\it Weak Closures of Ergodic Actions},
 Math. Notes, 101:2 (2017), 277-283.
 
 \bibitem{Mac}
  MacLane, S.
 {\it An algebra of additive relations.}
  Proc. Nat. Acad. Sci. USA,
  47 (1961) 1043-1051.

\bibitem{Mich}
Michael E.
{\it Topologies on spaces of subsets.}
Trans. Amer. Math. Soc. 71 (1951), 152--182.

\bibitem{Mack}
Mackey G. W.,
{\it Borel structure in groups and their duals.}
Trans. Amer. Math. Soc. 85 (1957), 134-165.


\bibitem{Nel}
 Nelson E., {\it The free Markoff field}, J. Funct. Anal. 12 (1973) 211-227.

\bibitem{Ner-bist}
 Neretin Yu. A.,
{\it Categories of bistochastic measures, and representations of some infinite-dimensional groups},
Russian Acad. Sci. Sb. Math., 75:1 (1993), 197-219.

\bibitem{Ner-book}
Neretin Yu.A., {\it Categories of Symmetries and Infinite-Dimensional Groups}, Oxford University Press, New York, 1996.

\bibitem{Ner-matching}
 Neretin Yu. A., {\it Spreading maps {\rm (}polymorphisms{\rm)}, symmetries of Poisson processes, and matching summation.}
 J. Math. Sci. (N. Y.), 126:2 (2005), 1077-1094.
 
 \bibitem{Ner-cl}
  Neretin, Yu.
 {\it Lectures on Gaussian integral operators and classical groups.}
   European Mathematical Society (EMS), Z\"urich, 2011
 
 \bibitem{Ner-gauss}
 Neretin, Yu. {\it Symmetries of Gaussian measures and operator colligations.} J. Funct. Anal. 263 (2012), no. 3, 782-802.
 
 \bibitem{Ner-finite}
  Neretin Yu. A. {\it The space $L^2$ on semi-infinite Grassmannian over finite field.} Adv. Math. 250 (2014), 320-350. 
 
 \bibitem{Ner-Haar}
  Neretin  Yu. A.,
{\it Topological groups and invariant measures.} Preprint,  arXiv:1510.03082. 
 
 
 \bibitem{Ner-Whi}
 Neretin  Yu. A., {\it Wishart--Pickrell distributions and closures of group actions}, J. Math. Sci. (N. Y.), 224:2 (2017), 328-334. 
 
  \bibitem{Ner-Wei}
 Neretin  Yu. A.,
 {\it On the Weil representation of infinite-dimensional symplectic group over a finite field.}
 Preprint
 \newline
  {\tt https://arxiv.org/abs/1703.07238}.
 
 \bibitem{Ner-GL}
  Neretin  Yu. A.,
 {\it Groups $\GL(\infty)$ over  finite fields and multiplications of double cosets.}
 Preprint {\tt https://arxiv.org/abs/2002.09969}
 
 \bibitem{Nev}
  Neveu, J.
 {\it  Mathematical foundations of the calculus of probability.}
  San Francisco, Holden-Day, Inc., 1965
 
 \bibitem{Ols-GB}
 Olshanski  G.I., {\it Unitary representations of infinite dimensional pairs (G,K) and the formalism of R. Howe}, in:
 {\it Representation of Lie Groups and Related Topics} (eds.  Vershik A. M.  and Zhelobenko D. P.),
 Gordon \& Breach 7 (1990) 269-463.
 
 
 
 
 \bibitem{Ols-semi}
 Olshanski G.I., {\it On semigroups related to infinite-dimensional groups,} in: A.A. Kirillov (Ed.),
 {\it Topics in Representation Theory}, Amer. Math. Soc., Providence, RI, 1991, pp. 67--101.
 
 
  
 
 \bibitem{Ryzh}
Ryzhikov  V. V.  {\it Weak Closure of Infinite Actions of Rank 1, Joinings, and Spectrum}, Math. Notes, 106:6 (2019), 957-965.  



\bibitem{SV}
Schmidt K., Vershik A. {\it Algebraic polymorphisms.} Ergodic Theory Dynam. Systems 28 (2008), no. 2, 633-642.



\bibitem{ShG}
Shilov G. E., Gurevich B. L. {\it Integral, measure and derivative: A unified approach.}
 Prentice-Hall, Inc.,
Englewood Cliffs, N.J. 1966 



\bibitem{Sch}
Schubert H.
{\it Categories.} 
Berlin-Heidelberg-New York: Springer-Verlag, 1972. 




\bibitem{Sol}
Solecki S.
{\it Closed subgroups generated by generic measure automorphisms.}
Ergodic Theory Dynam. Systems 34 (2014), no. 3, 1011-1017.



\bibitem{Tsa}
Tsankov, T. {\it Unitary representations of oligomorphic groups.} Geom. Funct. Anal. 22 (2012), no. 2, 528--555. 

 

\bibitem{Ver}
 Vershik A. M., {\it Many-valued measure-preserving mappings (polymorphisms) and Markovian operators},
J. Soviet Math., 23:3 (1983), 2243-2266.



\bibitem{Wei}
Weil, A. {\it L'int\'gration dans les groupes topologiques et ses applications.} (French)  Hermann et Cie., Paris, 1940.

 

 
 \bibitem{Wolf}
  Wolf J. A.,
 {\it Elliptic spaces in Grassmann manifolds.}
 Ill. J. Math., 7 (1963)
 447-462.
 
 
 
 \bibitem{Zhe}
 Zhelobenko, D. P.
 {\it  Principal structures and methods of representation theory.}
  Providence, RI,
  American Mathematical Society (AMS),
  2006
 
\end{thebibliography}
\end{document}